\documentclass{elsarticle}
\usepackage{graphicx}
\usepackage{amsmath,amssymb}
\usepackage{txfonts}
\usepackage{subfigure}
\usepackage{pslatex}
\usepackage{multirow}
\usepackage{colortbl}
\usepackage{wasysym}
\usepackage{verbatim}
\newtheorem{theorem}{Theorem}[section]
\newtheorem{definition}[theorem]{Definition}

\newenvironment{defn*}{\begin{definition}}{\end{definition}}

\numberwithin{equation}{section}

\begin{document}

\begin{frontmatter}
\title{Energy Stable Discontinuous Galerkin Finite Element Method for the Allen--Cahn Equation}

\author[bulent]{B\"{u}lent Karas\"{o}zen}
\author[ayse]{Ay\c{s}e Sar{\i}ayd{\i}n Filibelio\u{g}lu\corref{cor}}
\author[murat]{Murat Uzunca}

\address[bulent]{Department of Mathematics and Institute of Applied
  Mathematics, Middle East Technical University, 06800 Ankara,
  Turkey}
\address[ayse]{Institute of Applied Mathematics,
        Middle East Technical University, 06800 Ankara, Turkey}
\address[murat]{Department of Mathematics, At{\i}l{\i}m University, 06836 Ankara, Turkey and  Institute of Applied Mathematics,
        Middle East Technical University, 06800 Ankara,
  Turkey}
\cortext[cor]{Corresponding author.}

 \begin{abstract}
Allen--Cahn equation with constant and degenerate mobility, and with polynomial and logarithmic energy functionals is discretized using symmetric interior penalty discontinuous Galerkin (SIPG) finite elements in space.  We show that the energy stable average vector field (AVF) method as the time integrator for gradient systems like the Allen-Cahn equation satisfies the energy decreasing property for the fully discrete scheme. The numerical results for one and two dimensional Allen-Cahn equation with periodic boundary condition, using adaptive time stepping,  reveal that  the discrete energy decreases monotonically, the phase separation and metastability phenomena can be observed and the ripening time is detected correctly.
\bigskip

 \end{abstract}
 \end{frontmatter}

\section{Introduction}
In gradient flows, the energy of the system  decreases along the solutions as fast as possible. A typical example is the Allen--Cahn   equation modeling  the reaction-diffusion process in material sciences
\begin{equation} \label{acgradient}
u_t = -\mu(u)\ \frac{\delta \mathcal{E}(u)}{\delta u}
\end{equation}
with minimizing the Ginzburg--Landau energy functional
\begin{equation}\label{energy}
\mathcal{E}(u)=\int_{\Omega} \left( \frac{\epsilon^2}{2}|\nabla u|^2 + F(u)  \right)dx.
\end{equation}
The term  $ \frac{\delta \mathcal{E}(u)}{\delta u}$ in \eqref{acgradient} denotes the variational derivative of \eqref{energy} in the $L^2$ norm in the domain $\Omega \subset \mathbb{R}^d$ ($d=1,2$).
The Allen--Cahn   equation 
\begin{equation}\label{allencahn}
u_t=\mu(u)(\epsilon^2\Delta u -f(u)), \quad (x,t) \in \Omega  \times (0,T]
\end{equation}
was first introduced by Allen and Cahn \cite{allen79amt} to describe the motion of anti-phase boundaries of a binary alloy at a fixed temperatures. In the last thirty years, Allen-Chan equation has been widely used in many complicated moving interface problems in material science, fluid dynamics, image analysis and mean curvature flow. In \eqref{allencahn}, the unknown $u$  denotes the concentration one of the species of the alloy, known as the phase state between materials. The parameter $\epsilon$ is the interaction length, capturing the dominating effect of the reaction kinetics and represents the effective diffusivity and $\mu(u)$ is the non--negative mobility function that describes the physics of phase separation. The nonlinear term $f(u)=F^{'}(u)$ is the derivative of a free energy functional $F(u)$. Two types of free energy functional $F(u)$ are considered in the literature.
The first one is the non-convex logarithmic free energy  \cite{barrett99fea,barrett00fea,shen14omp}
\begin{equation}\label{func2}
F(u)=\frac{\theta}{2}[(1+u)\ln(1+u) + (1-u)\ln(1-u)] - \frac{\theta_c}{2}u^2
\end{equation}
with $ 0 < \theta \le \theta_c$, where $\theta_c$ is the transition temperature. For temperatures $\theta$ close to $\theta_c$, the logarithmic free energy  (\ref{func2}) is usually approximated by the convex quartic double-well potential  \cite{feng13nsi,willoughby11hor}
\begin{equation}\label{func1}
F(u)=\frac{1}{4}(1-u^2)^2.
\end{equation}
In the case of the quartic double-well potential \eqref{func1}, $f(u)=u^3-u$  represents the bi-stable non-linearity. For the logarithmic free energy (\ref{func2}) it takes the form $f(u)=\frac{\theta}{2}\ln\left(\frac{1+u}{1-u}\right) - \theta_c u$.
The logarithmic free energy and degenerate mobility are often used for the Cahn-Hilliard equation \cite{barrett99fea,barrett00fea,Guo14esd}. Common choices for the degenerate mobility are  $\mu (u)=\beta (1-u^2)$ or  $\mu (u)=\beta u(1-u)$ with constant   $\beta$.
In the literature, the Allen-Cahn equation is investigated using the constant mobility and the quartic double-well potential.  Allen-Cahn equation with degenerate mobility and logarithmic free energy was introduced first time in \cite{shen14omp}.

The mobility function $\mu(u)$, and both the free energy functionals (\ref{func2}) and (\ref{func1}) together with their derivatives  are Lipschitz continuous for $u_1, u_2 \in \mathbb{R}$ with the constraints $|u_{1,2}|\le 1$  \cite{zee11goe}:
\begin{eqnarray}
\left|  \mu(u_1) - \mu(u_2)\right|\leq L_{b} \left| u_1 - u_2\right|, \nonumber \\
\left|  f(u_1) - f(u_2)\right|\leq L_{f} \left| u_1 - u_2\right|,  \\
\left|  f'(u_1) - f'(u_2)\right|\leq L_{f'} \left| u_1 - u_2\right|, \nonumber
\end{eqnarray}
where $L_{\mu}, L_{f}, L_{f'}\geq 0$ stand for the related Lipschitz constants.

Energy decrease property of the Allen-Cahn equation is obtained by  taking the $L^2$-inner product of (\ref{allencahn}) with $(-\epsilon^2\Delta u + f(u))$
\begin{equation}\label{energydecrrease}
\mathcal{E}(u(t_{n})) < \mathcal{E}(u(t_{m})), \quad \forall t_{n} > t_{m}.
\end{equation}

The presence of the small inter-facial length $\epsilon$, different time scales of the phase separation  and coarsening, the non-linearity are the main challenges in the numerical solution of the Allen-Cahn equation. For space discretization, well known methods like  finite-difference, spectral elements \cite{christlieb14has}, continuous finite element \cite{liu13ssi} and local discontinuous Galerkin (LDG) methods \cite{guo-nse,feng14aip}  are used. Several energy stable integrators are developed to preserve the energy decreasing property of the Allen-Cahn equation with constant mobility. For small values of the diffusion parameter $\epsilon$, semi-discretization in space leads to stiff systems. Because the explicit methods are not suitable for stiff systems and the fully implicit systems require solution of non-linear equations at each time step, implicit-explicit (IMEX) methods \cite{shen10nac} and  parametrized energy stable semi-implicit schemes are developed \cite{shen10nac,feng13scn,feng13nsi}.

In this work, we use the symmetric interior penalty discontinuous Galerkin finite elements (SIPG) for space discretization \cite{arnold82ipf,riviere08dgm} and the energy stable average vector field (AVF) integrator for time discretization. In contrast to the continuous finite elements, the discontinuous finite elements use piecewise polynomials that are fully discontinuous at the interfaces. In this way, the SIPG approximation allows to capture the sharp gradients or singularities locally. It is important to design efficient and accurate numerical schemes that are energy stable  and robust for small $\epsilon$. Among the energy stable implicit methods, the best known is the  first order implicit Euler method which is strongly energy decreasing, i.e. the discrete energy decreases without any restriction for the step size $\Delta t$ for very stiff gradient systems for very small $\epsilon$ \cite{hairer10epv}. The only second order implicit energy stable method for the gradient systems is the average vector field (AVF) method  \cite{celledoni12ped,hairer10epv}.
The mid-point method coincides with the AVF method for quadratics nonlinearities. For gradient systems like the Allen-Cahn equation involving higher order polynomial or general nonlinear terms, the mid-point method is not energy stable   \cite{hairer10epv}.
Higher order energy decreasing methods are the discontinuous Galerkin-Petrov in time methods (with different trial and test functions) \cite{schwieweck10sdg} and Gauss Radau IIA Runge-Kutta collocation methods \cite{hairer96sod}. But they require coupled systems of equations at each time step and are computationally cost.

The rest of the paper is organized as follows.  In Section \ref{dgfem}, we give the SIPG semi-discretization of the Allen-Cahn equation with degenerate mobility for periodic boundary conditions.  Section \ref{avf} is devoted to the time discretization with the AVF method, where the solution of the system of non-linear equations are described. The non-linear energy stability of the fully discrete scheme is given in Section \ref{stability}. We present in Section \ref{numerical} several numerical examples to demonstrate the performance of the SIPG discretization coupled with AVF method for the Allen-Cahn equation using a time-adaptive algorithm. The paper ends with some conclusions.


\section{Symmetric interior penalty Galerkin discretization}\label{dgfem}
In this section, we briefly describe the symmetric interior penalty Galerkin (SIPG) discretization of the Allen-Cahn equation (\ref{allencahn}) equipped with periodic boundary conditions on a 2D domain $\Omega\subset\mathbb{R}^2$. The classical (continuous) weak formulation (semi--discrete) of the Allen-Cahn equation (\ref{allencahn}) reads as: find $u(t)\in H_{per}^{1}(\Omega)$ such that
\begin{equation} \label{2}
\begin{aligned}
(u_t, v)_{\Omega} + a(u,v)+ (\tilde{f}(u),v)_{\Omega} &= 0\; , \qquad \qquad \forall v\in H_{per}^{1}(\Omega), \quad t\in (0,T],\\
(u(0),v)_{\Omega} &= (u_0,v)_{\Omega} \; , \; \quad \forall v\in H_{per}^{1}(\Omega),
\end{aligned}
\end{equation}
where $(\cdot , \cdot)_{\Omega}$ denotes the usual $L^2$-inner product over the domain $\Omega$, the bilinear form $a(u,v)=\epsilon^2\mu (u)(\nabla u,\nabla v)_{\Omega}$, the nonlinear term $\tilde{f}(u)=\mu (u)f(u)$ and $H_{per}^{1}(\Omega)$ is the restricted finite element space given by
$$
H_{per}^{1}(\Omega)=\{ v\in H^{1}(\Omega):\; v|_{\Gamma_1}=v|_{\Gamma_2}\},
$$
where $\Gamma_1\subset\partial\Omega$ and $\Gamma_2\subset\partial\Omega$ are the subsets of the domain boundary related to the parts of periodic boundary pairs satisfying $\Gamma_1\cap\Gamma_2 = \emptyset$ and $\Gamma_1\cup\Gamma_2 =\partial\Omega$. In the case of degenerate mobility, we take the mobility function computed from the previous time step (explicitly) as for the Cahn-Hilliard equation  \cite{barrett99fea,barrett00fea,Guo14esd}. In the classical (continuous) finite elements, the time-dependent approximation to the system (\ref{2}) belongs to a finite dimensional conforming subspace $V_h\subset H_{per}^{1}(\Omega)$. In contrast to the continuous finite elements, discontinuous Galerkin (DG) methods use non-conforming spaces, i.e. $V_h\not\subset H_{per}^{1}(\Omega)$ and there is no need for a restriction on the FE space.

The DG discretization given in this article is based on the  SIPG method \cite{riviere08dgm} applied to the diffusion part of the problem equipped with periodic boundary conditions \cite{vemaganti07dgm}. Let $\left\{\mathcal{T}_{h}\right\}_{h}$ be a family of shape regular meshes such that $\bar{\Omega}=\cup_{K \in \mathcal{T}_h}\bar{K}$, $K_i\cap K_j=\emptyset$ for $K_i$, $K_j \in \mathcal{T}_{h}$,  $i\neq j$. The diameter of an element $K$ and the length of an edge $E$ are denoted by $h_K$ and $h_E$, respectively. Set the test and trial space 
\begin{equation}
V_h = \left\{u \in L^{2}(\Omega) : u|_{K} \in \mathbb{P}^{q}(K) \quad \forall K \in  \mathcal{T}_{h} \right\},
\end{equation}
where $\mathbb{P}^{q}(K)$ denotes the set of all polynomials on $K \in  \mathcal{T}_{h}$ of degree at most $q$. We note that the trial space and the space of test functions are chosen to be the same without, in contrast to continuous finite elements, any restriction on the boundary.

We split the set of all edges $E_h$ into the set $E^{0}_{h}$ of interior edges and the set $E^{per}_{h}$ of periodic boundary edge--pairs. An individual element of the set $E^{per}_{h}$ is of the form $\omega =\{ E_l, E_m\}$ where $E_l \subset \partial K_{l}\cap \partial\Omega$, and $E_m \subset \partial K_{m}\cap \partial\Omega$ is the corresponding periodic edge-pair of $E_l$  with $l>m$, and we associate with each $\omega$ a common normal vector $\mathbf{n}$ that is outward unit normal to $E_l\subset \partial K_{l}\cap \partial\Omega$. Let the edge $E$ be a common interior edge for two elements $K$ and $K^{e}$. For a piecewise continuous scalar function $u$, because of the discontinuity on the interfaces, there are two traces of $u$ along $E$, denoted by $u|_{E}$ from inside $K$ and $u^{e}|_{E}$ from inside $K^{e}$. Then, define the jump and average of $u$ across the edge $E$ as:
$$
[u]=u|_{E}\mathbf{n}_{K} + u^{e}|_{E}\mathbf{n}_{K^{e}}, \quad \{u\}=\frac{1}{2}(u|_{E} + u^{e}|_{E}).
$$
where $\mathbf{n}_{K}$ and $\mathbf{n}_{K^e}$ denote the outward unit normal vector to the boundary of the elements $K$ and $K^e$ on the edge $E$, respectively. Similarly, for a piecewise continuous vector valued function $\nabla u$, the jump and average across an edge $E$ are given by
$$
[\nabla u]=\nabla u|_{E} \cdot \mathbf{n}_{K} + \nabla u^{e}|_{E} \cdot \mathbf{n}_{K^{e}}, \quad \{\nabla u\}=\frac{1}{2}(\nabla u|_{E} + \nabla u^{e}|_{E}).
$$
In case of the boundary edges, the periodic boundary edges are treated as interior edges, in other words,  as unknown with appropriate definitions of the so--called jump and average terms. Then, for each $\omega=\{ E_l, E_m\}\in E^{per}_{h}$, we define the jump and average operators as
$$
[u]_{\omega}=u|_{E_l}\mathbf n - u|_{E_m}\mathbf n, \quad \{ u\}_{\omega}=\frac{1}{2}(u|_{E_l} + u|_{E_m)}.
$$

Following the definitions above, the SIPG semi-discretized system of the Allen-Cahn equation \eqref{allencahn} reads as: set $u_h(0)\in V_h$ be the projection (orthogonal $L^2$-projection) of the initial condition $u_0$ onto $V_h$, find $u_h(t)\in V_h$ such that
\begin{equation}\label{discreteallencahn}
(\partial_{t}u_h, \upsilon_h)_{\Omega} + a_{h}(\epsilon^2\mu (u_h);u_{h},\upsilon_{h}) + (\tilde{f}(u_{h}),\upsilon_{h})_{\Omega}=0,\: \quad \forall \upsilon_{h} \in V_h, \: t \in (0,T],
\end{equation}
where $a_{h}(\kappa; u,\upsilon)= \tilde{a}_{h}(\kappa; u,\upsilon) + J_{h}^{\partial}(\kappa;u,\upsilon)$ is the bilinear form with
\begin{align*}
\tilde{a}_{h}(\kappa; u,\upsilon ) &= \sum_{K \in \mathcal{T}_h}\int_K \kappa \nabla u \cdot\nabla \upsilon - \sum_{E\in E^{0}_{h}}\int_E\left\{\kappa\nabla u\right\}\cdot[\upsilon]ds\\
 & \qquad - \sum_{E\in E^{0}_{h}}\int_E\left\{\kappa\nabla \upsilon \right\}\cdot[u]+\sum_{E\in E^{0}_{h}} \frac{\sigma \kappa}{h_{E}}\int_E [u]\cdot[\upsilon]ds,\\
J_{h}^{\partial}(\kappa;u,\upsilon) &= - \sum_{\omega\in E^{per}_{h}}\int_{\omega}\left\{\kappa \nabla u\right\}_{\omega}\cdot[\upsilon]_{\omega}ds - \sum_{\omega\in E^{per}_{h}}\int_{\omega}\left\{\kappa \nabla \upsilon \right\}_{\omega}\cdot[u]_{\omega} \\
& \qquad +\sum_{\omega\in E^{per}_{h}} \frac{\sigma \kappa}{h_{E}}\int_{\omega} [u]_{\omega}\cdot[\upsilon]_{\omega}ds
\end{align*}

The parameter $\sigma$ in the above formulation is called the penalty parameter and it should be sufficiently large to ensure the stability of the DG discretization with a lower bound depending only on the polynomial degree $q$, for instance, for 1D models usually $\sigma=2.5(q+1)^2$ is taken \cite{riviere08dgm}.

\section{Time discretization by average vector field method}\label{avf}

The semi-discrete system \eqref{discreteallencahn} of the Allen-Cahn equation can also be regarded as a gradient system  of the form $\dot{y} = -\nabla U(y)$  evolving into a state of minimal energy, characterized by the monotonically energy decreasing property of the potential
$$
U(y(t)) \le U(y(s)), \quad \hbox{ for } t > s.
$$
In the numerical approximation of the gradient systems, it is desirable to preserve the energy decreasing property monotonically
$$
U(y(t_n)) \le U(y(t_{n-1})), \quad \hbox{ for } n=1,2\ldots
$$

The average vector field (AVF) method \cite{celledoni12ped,hairer10epv}
$$
y_n = y_{n-1} - \Delta t \int _{0}^{1} \nabla U (\tau y_n+( 1-\tau )y_{n-1})d \tau
$$
possesses the energy decreasing property without restriction to  step size $\Delta t$.  It represents a modification of the implicit mid-point rule and  for quadratic potentials  $U(y)$, the AVF reduces to the mid-point rule.

The AVF method is  also equivalent to the Petrov-Galerkin discontinuous Galerkin in time, when the trial functions are piecewise linear and the test functions are piecewise constant, which are given by
\begin{equation} \label{dgavf}
y_{n}=y_{n-1}-\int_{t_{n-1}}^{t_n}\nabla U(\Delta t^{-1}(t-t_{n-1})y_n+\Delta t^{-1}(t_n-t)y_{n-1})dt.
\end{equation}
With the time parametrization $t(\tau ) = t_{n-1} + (t_n-t_{n-1})\tau $
and using the change of variables formulation
$$
\int_{t_{n-1}}^{t_n}g(t)dt = \int_0^1g(t(\tau ))\frac{dt(\tau )}{d\tau }d\tau,
$$
we obtain for the integral term in (\ref{dgavf})
\begin{eqnarray*}
y_{n} &=& y_{n-1}-\int_{t_{n-1}}^{t_n}\nabla U(\Delta t^{-1}(t-t_{n-1})y_n+\Delta t^{-1}(t_n-t)y_{n-1})dt\\
&=& y_{n-1}-\Delta t\int_0^1 \nabla U(\tau y_n+( 1-\tau )y_{n-1})d\tau ,
\end{eqnarray*}
which is the AVF method on the interval $[t_{n-1},t_n]$.

The SIPG semi-discretized system \eqref{discreteallencahn} of the Allen-Cahn equation has the time-dependent solution of the form
\begin{equation}\label{4}
u_h(t)=\sum^{N}_{m=1}\sum^{n_q}_{j=1}\xi^{m}_{j}(t) \varphi^{m}_{j},
\end{equation}
where $\varphi^{m}_{j}$ and $\xi^{m}_{j}$, $j=1, \ldots, n_q$, $m=1, \ldots, N$, are the basis functions spanning the space $V_h$ and the unknown coefficients, respectively. The number $n_q$ denotes the local dimension of each DG element (interval in 1D, triangle in 2D) with  $n_q=q+1$ for 1D problems and $n_q=\frac{(q+1)(q+2)}{2}$ for 2D problems, and $N$ is the number of DG elements. Substituting \eqref{4} into (\ref{discreteallencahn}) and choosing $\upsilon=\varphi^{k}_{i}, \: i=1, \ldots, n_q$, $k=1, \ldots, N$, we obtain the following semi-linear system of ordinary differential equations as a gradient system
\begin{equation}\label{7}
M\xi_{t}=-\nabla U(\xi)=- A\xi - r(\xi )
\end{equation}
for the unknown coefficient vector $\xi =(\xi_{1}^{1}, \ldots , \xi_{n_q}^{1},\xi_{1}^{2},\ldots , \xi_{1}^{N}, \ldots , \xi_{n_q}^{N})^T$, where $M$ is the mass matrix, $M_{ij}=(\varphi^{j}, \varphi^{i})_{\Omega}$, $1\leq i,j \leq n_q\times N$, $A$ is the stiffness matrix, $A_{ij}=a_h(\kappa;\varphi^{j}, \varphi^{i})$, $1\leq i,j \leq n_q\times N$, $r$ is the non--linear vector of unknown coefficient vector $\xi$ with the entries $r_{i}(\xi)=(\tilde{f}(u_h),\varphi^{i})_{\Omega}$, $1\leq i \leq n_q\times N$.

We consider the uniform partition $0=t_0<t_1<\ldots < t_J=T$ of the time interval $[0,T]$ with the uniform time step-size $\Delta t=t_{k}-t_{k-1}$, $k=1,2,\ldots , J$. We set $\xi_{n}\approx \xi (t_n)$ as the approximate solution at the time instance $t=t_n$, $n=0,1,\ldots ,J$. For $t=0$, let $u_h(0)\in V_h$ be the projection (orthogonal $L^2$-projection) of the initial condition $u_0$ onto $V_h$, and let $\xi_0$ be the corresponding coefficient vector  satisfying (\ref{4}). Then, the average  vector field method applied to the gradient system (\ref{7}) reads as: for $n=0,1,\ldots , J-1$, solve
\begin{eqnarray*}
\frac{M\xi_{n+1}-M\xi_{n}}{\Delta t}&=&  -\int^{1}_{0} \nabla U( \tau \xi_{n+1} + (1-\tau)\xi_n )d\tau\\
M\xi_{n+1}&=&M\xi_{n} - \Delta t\underbrace{\int^{1}_{0} \left[ A ( \tau \xi_{n+1} + (1-\tau)\xi_n ) \right] d\tau}_{linear} - \Delta t\underbrace{\int^{1}_{0} r( \tau \xi_{n+1} + (1-\tau)\xi_n ) d\tau}_{non-linear}.\\
\end{eqnarray*}
After a simple calculation for the linear part, we get
\begin{equation}\label{8}
M\xi_{n+1}=M\xi_{n} - \frac{\Delta t}{2}(A\xi_n + A\xi_{n+1}) - \Delta t\int^{1}_{0} r( \tau \xi_{n+1} + (1-\tau)\xi_n ) d\tau ,
\end{equation}
which is the fully discretized system that we will solve for $\xi_{n+1}$. We solve the non-linear system of equations \eqref{8} using Newton's method. From the algebraic point of view, Newton's method for (\ref{8}) corresponds to solving the residual equations

\begin{equation}\label{nonlin}
R(\xi_{n+1})= M\xi_{n+1} - M\xi_{n} +\frac{\Delta t}{2}(A\xi_n + A\xi_{n+1}) + \Delta t\int^{1}_{0} r( \tau \xi_{n+1} + (1-\tau)\xi_n ) d\tau =0.
\end{equation}

Starting with an initial guess $\xi_{n+1}^{(0)}$, the $k-{th}$ Newton iteration to solve the residual equation (\ref{nonlin}) for the unknown vector $\xi_{n+1}$ reads as
\begin{equation}\label{newt}
Js^{(k)} = - R(\xi_{n+1}^{(k)}), \qquad
\xi_{n+1}^{(k+1)} = \xi_{n+1}^{(k)} + s^{(k)} \; , \quad k=0,1,\ldots
\end{equation}
until a user defined tolerance is satisfied. In (\ref{newt}), $J$ stands for the Jacobian matrix of $R(\xi_{n+1})$, whose entries are the partial derivatives
$$
J_{ij}=\frac{\partial R_i}{\partial (\xi_{n+1})_j} \; , \qquad i,j=1,2,\ldots , n_q\times N
$$
at the current iterate $\xi_{n+1}^{(k)}$. It is easy to differentiate the linear terms in (\ref{nonlin})
$$
\frac{\partial }{\partial (\xi_{n+1})_j}\left( M\xi_{n+1} - M\xi_{n} + \frac{\Delta t}{2}(A\xi_n + A\xi_{n+1}) \right)_i = M_{ij} + \frac{\Delta t }{2}A_{ij}.
$$
To differentiate the non-linear term in (\ref{nonlin}), we apply the chain rule
\begin{eqnarray*}
\frac{\partial }{\partial (\xi_{n+1})_j}\Delta t\int^{1}_{0} r_i( \hat{\xi} ) d\tau &=& \Delta t\int^{1}_{0} \tau \frac{\partial r_i( \hat{\xi} )}{\partial \hat{\xi}_j}  d\tau ,
\end{eqnarray*}
where we have set $\hat{\xi}=\tau \xi_{n+1} + (1-\tau)\xi_n$, and using the expansion $\hat{u}_h=\sum^{n_q\times N}_{k=1} \hat{\xi}_{k} \varphi^k$, ordered version of (\ref{4}), we get 
\begin{eqnarray}\label{jacobnonlin}
\frac{\partial r_i(\hat{\xi} )}{\partial \hat{\xi}_j}&=&\frac{\partial }{\partial \hat{\xi}_j}  (\tilde{f}(\hat{u}_h), \varphi^i)_{\Omega},\qquad i,j=1,2,\ldots , n_q\times N \nonumber \\
                                         &=& \frac{\partial }{\partial \hat{\xi}_j}  (\mu(u^{n}_{h})f(\hat{u}_h), \varphi^i)_{\Omega} \nonumber \\
                                         &=& \mu(u^{n}_{h})\int_{\Omega} f'\left( \sum^{n_q\times N}_{j=1}\hat{\xi}_{j} \varphi^j\right) \varphi^j\varphi^i dx.
\end{eqnarray}
We obtain finally the Jacobian matrix as

\begin{equation}\label{jacob}
J= M + \frac{\Delta t }{2}A + \Delta t  \int_0^1 \tau J_{r}( \tau \xi_{n+1}^{(k+1)} + (1-\tau)\xi_n ) d\tau ,
\end{equation}
where $J_{r}( \tau \xi_{n+1}^{(k+1)} + (1-\tau)\xi_n )$ is the differential matrix, whose entries are given in (\ref{jacobnonlin}) for $ \hat{\xi}=\tau \xi_{n+1}^{(k+1)} + (1-\tau)\xi_n$. At each Newton iteration, we approximate the integral term in (\ref{jacob}) using the fourth order Gaussian quadrature rule, by which, in the case of polynomial nonlinear terms arising from the double-well potential \eqref{func2}, the integrals are evaluated exactly.

\section{Energy stability of the fully discrete scheme}\label{stability}

In this Section, we show the non-linear stability of the AVF method applied to the semi-discrete system \eqref{discreteallencahn}. The DG discretized  energy of the semi-discrete Allen--Cahn equation \eqref{discreteallencahn} at the time $t^n=n \Delta t$ is given as \cite{feng14aip}

\begin{equation}\label{discreteenergy}
\mathcal{E}^h_{DG}(u^{n}) = \frac{\epsilon^2}{2}\left\|\nabla u^{n}\right\|^2_{L^2(\tau_h)} + (F(u^n),1)_{\Omega} + \sum_{E\in E^{0}_{h}} \left( -(\{\epsilon^2\partial _n u^{n}\}, [u^{n}])_E + \frac{\sigma\epsilon^2}{2h_E}([u^{n}],[u^{n}])_E \right).
\end{equation}

Time discretization of semi-discrete system \eqref{discreteallencahn} by the AVF methods, leads to

\begin{eqnarray*}
\frac{1}{\Delta t }(u^{n+1}_{h} - u^{n}_{h}, \upsilon_{h})_{\Omega} & + &\mu(u^{n}_{h})\frac{1}{2}a_h(\epsilon^2;u^{n+1}_{h}+u^{n}_{h},\upsilon_{h})\\
& + & \mu(u^{n}_{h})\int_0^1(f(\tau u^{n+1}_{h} + (1-\tau)u^{n}_{h}), \upsilon_{h})_{\Omega} d\tau = 0, \quad \forall \upsilon_{h} \in V_h.
\end{eqnarray*}

Taking $\upsilon_{h}=u^{n+1}_{h}-u^{n}_{h}$, we obtain

\begin{eqnarray*}
\frac{1}{\Delta t}(u^{n+1}_{h} - u^{n}_{h}, u^{n+1}_{h}-u^{n}_{h})_{\Omega} & + & \mu(u^{n}_{h})\frac{1}{2}a_h(\epsilon^2;u^{n+1}_{h}+u^{n}_{h},u^{n+1}_{h}-u^{n}_{h}) \\
 & + & \mu(u^{n}_{h})\int_0^1(f(\tau u^{n+1}_{h} + (1-\tau)u^{n}_{h}), u^{n+1}_{h}-u^{n}_{h})_{\Omega} d\tau =0.
\end{eqnarray*}

By using the identity $(a+b, a-b)_{\Omega}= (a^2 - b^2 ,1)_{\Omega}$ and the bilinearity of $a_h$, we get

\begin{align} \label{eq2avff}
\frac{1}{ \Delta t}(u^{n+1}_{h} - u^{n}_{h}, u^{n+1}_{h}-u^{n}_{h})_{\Omega} + \mu(u^{n}_{h})\int_0^1(f(\tau u^{n+1}_{h} + (1-\tau)u^{n}_{h}), u^{n+1}_{h}-u^{n}_{h})_{\Omega} d\tau \nonumber \\
\qquad + \mu(u^{n}_{h})\frac{1}{2}a_h(\epsilon^2;u^{n+1}_{h},u^{n+1}_{h}) - \mu(u^{n}_{h})\frac{1}{2}a_h(\epsilon^2;u^{n}_{h},u^{n}_{h}) =0.
\end{align}
Taylor expansions of $F$ around $u^{n}_{h}$ and $u^{n+1}_{h}$ leads to
\begin{align*}
F(u^{n}_{h}) &\approx F(\tau u^{n+1}_{h} + (1-\tau)u^{n}_{h}) - f(\tau u^{n+1}_{h} + (1-\tau)u^{n}_{h})(\tau (u^{n+1}_{h} - u^{n}_{h}))  \\
F(u^{n+1}_{h})&\approx  F(\tau u^{n+1}_{h} + (1-\tau)u^{n}_{h}) + f(\tau u^{n+1}_{h} + (1-\tau)u^{n}_{h})(1-\tau) (u^{n+1}_{h} - u^{n}_{h})).
\end{align*}

Subtracting $F(u^{n}_{h})$ from $F(u^{n+1}_{h})$ and ignoring higher order terms including the the derivatives of $f$, we obtain
\begin{align}
F(u^{n+1}_{h})-F(u^{n}_{h}) &\approx  f(\tau u^{n+1}_{h} + (1-\tau)u^{n}_{h}) (u^{n+1}_{h} - u^{n}_{h}) \nonumber \\
(F(u^{n+1}_{h}),1)_{\Omega}-(F(u^{n}_{h}),1)_{\Omega} &\approx  ( f(\tau u^{n+1}_{h} + (1-\tau)u^{n}_{h}), u^{n+1}_{h} - u^{n}_{h})_{\Omega} \nonumber \\
\int_0^1((F(u^{n+1}_{h}),1)_{\Omega}-(F(u^{n_{h}}),1)_{\Omega})d\tau &\approx \int_0^1( f(\tau u^{n+1}_{h} + (1-\tau)u^{n}_{h}), u^{n+1}_{h} - u^{n}_{h})_{\Omega}d\tau \nonumber \\
(F(u^{n+1}_{h}),1)_{\Omega}-(F(u^{n}_{h}),1)_{\Omega} &\approx \int_0^1( f(\tau u^{n+1}_{h} + (1-\tau)u^{n}), u^{n+1}_{h} - u^{n}_{h})_{\Omega}d\tau .\label{tayavf}
\end{align}
We note  that the bilinear form $a_h(\epsilon^2; u^{n+1},u^{n+1})$ satisfies
\begin{eqnarray}
a_h(\epsilon^2; u^{n+1},u^{n+1}) &=&  \epsilon^2 \| \nabla u^{n+1}\|_{L^2(\Omega)}^2 -2\sum_{E \in E^{0}_{h}} \int_E \{\epsilon^2  \nabla u^{n+1}\}[ u^{n+1}]ds \nonumber \\
&  &  + \sum_{E \in E^{0}_{h}} \frac{\sigma \epsilon^2}{h_E} \|[u^{n+1}]\|_{L^2(E)}^2 \ge 0. \label{posa}
\end{eqnarray}
Since all the terms in \eqref{posa} are non-negative (see \cite[Sec. 2.7.1]{riviere08dgm} for positivity of edge integral term), we have $a_h(\epsilon^2; u^{n+1},u^{n+1})\geq 0$ and similarly $a_h(\epsilon^2; u^{n},u^{n})\geq 0$. Using these identities, positivity of the mobility $\mu(u^{n}_{h})$, and substituting \eqref{tayavf} into \eqref{eq2avff}, we obtain
\begin{eqnarray*}
0 \geq- \frac{1}{\Delta t} \left\|u^{n+1}_{h}-u^{n}_{h}\right\|_{L^2(\Omega)} &\approx&  \mu(u^{n}_{h})\left((F(u^{n+1}_{h}), 1)_{\Omega} +\frac{1}{2}a_h(\epsilon^2;u^{n+1}_{h},u^{n+1}_{h}) \right) \\
&& - \mu(u^{n}_{h})\left( (F(u^{n}_{h}),1)_{\Omega} + \frac{1}{2}a_h(\epsilon^2;u^{n}_{h},u^{n}_{h}) \right)\\
0 \geq- \frac{1}{\Delta t \mu(u^{n}_{h})} \left\|u^{n+1}_{h}-u^{n}_{h}\right\|_{L^2(\Omega)} &\approx&  (F(u^{n+1}_{h}), 1)_{\Omega} +\frac{1}{2}a_h(\epsilon^2;u^{n+1}_{h},u^{n+1}_{h})\\
& & - (F(u^{n}_{h}),1)_{\Omega} + \frac{1}{2}a_h(\epsilon^2;u^{n}_{h},u^{n}_{h})\\
&=& \mathcal{E}(u^{n+1}_{h}) - \mathcal{E}(u^{n}_{h}) ,
\end{eqnarray*}
which implies that $\mathcal{E}(u^{n+1}_{h}) \leq  \mathcal{E}(u^{n}_{h}). $
\section{Numerical results}\label{numerical}
In all numerical experiments, we have used linear polynomials to form the DG space. Only for the ripening time calculations, quadratic elements are used. We have considered in all examples periodic boundary conditions. Until forming the metastable state, the initial dynamics require small time steps as the transition layers are formed. During the metastable state, the dynamics changes not much, larger time steps are required. Therefore uniform time steps will be inefficient as shown in \cite{christlieb14has,willoughby11hor}.  We use adaptive time stepping to resolve the multiple time dynamics of the Allen--Cahn equation.

\subsection{Adaptive time stepping}
The transition layers of the Allen-Cahn equation move quickly from one unstable equilibrium to the other one by crossing the zero axis. The time where the solution takes its minimum value is named as the ripening time. The ripening time is computed for the Allen-Cahn equation  with periodic boundary conditions using adaptive time stepping  in  \cite{christlieb14has,willoughby11hor}.
For the construction of adaptive time grids, one needs a local error estimator. For local error estimation, two discrete solutions $u_\tau$, $\hat{u}_\tau$ of order $p+1$ and $p$ are computed  such that
$$u_\tau(\tau)=u(\tau) + O(\tau^{p+2}), \quad \hat{u}_\tau (\tau)=u(\tau) + O(\tau^{p+1})$$
with $\tau$ denoting the time step size $\Delta t$. Then
\begin{equation}\label{adaptime}
\hat{\epsilon}_\tau=\left\|u_\tau(\tau) - \hat{u}_\tau \right\|=C\tau^{p+1}
\end{equation}
is an estimator of the actual error $\hat{\epsilon}_\tau$ of $\hat{u}_\tau$ measured in an Euclidean norm  \cite{deuflhard12ans}. We search for an optimal step size $\tau^*$ for which $\hat{\epsilon}_{\tau^*}\leq \delta_{TOL}$, where $\delta_{TOL}$ denotes a user specified tolerance. By insertion of both $\tau$ and $\tau^*$ into (\ref{adaptime}), we arrive at the estimation formula
$$\tau^*=\sqrt[p+1]{\frac{\rho \delta_{TOL}}{\hat{\epsilon}_\tau}}\tau$$ with a safety factor $\rho\approx0.9$. If $\hat{\epsilon}_{\tau^*} \leq \delta_{TOL}$, then the presented step size $\tau^*$ is accepted and $\tau^*$ is used in the next step; otherwise the present step size is rejected and the current step is repeated with the step size $\tau^*$. In the successful case, the more accurate value $u_\tau(\tau)$ will be used to start the next step. For the  adaptive time stepping scheme, we choose backward Euler method and AVF method which are first and second order, respectively. We further set the initial time step size $\tau =0.05$ and $\delta_{TOL}=10^{-4}$.

\subsection{1D Allen--Cahn equation with constant mobility and double-well potential}\label{ex2}
\begin{figure}[htb!]
\centering
\includegraphics[height=0.4\textwidth, width=0.6\textwidth] {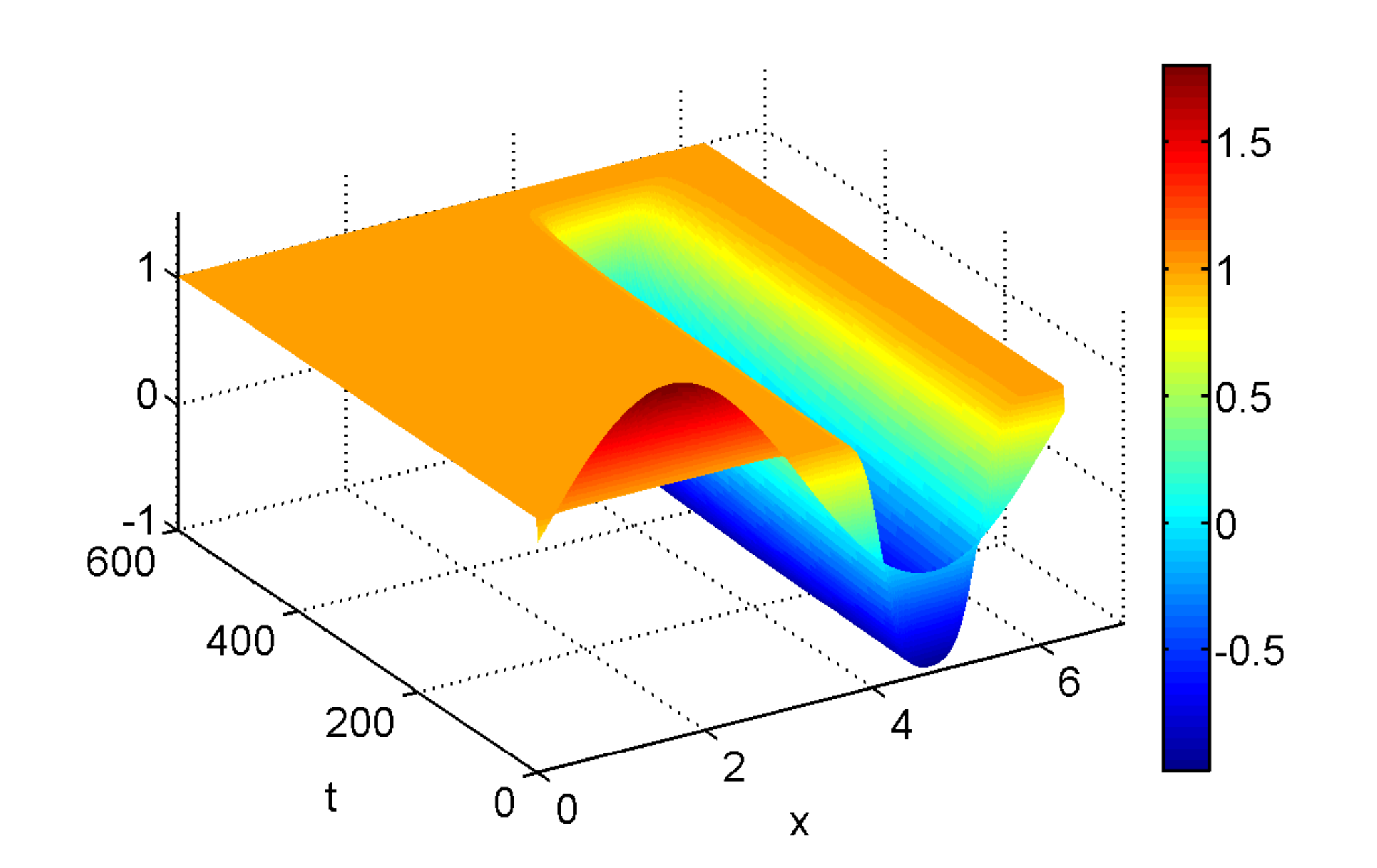}	
\caption{Example \ref{ex2}: The evolution of phase function.\label{1DPeriodic}}
\end{figure}

\begin{figure}[htb!]
	\centering
\includegraphics[height=0.4\textwidth, width=0.4\textwidth] {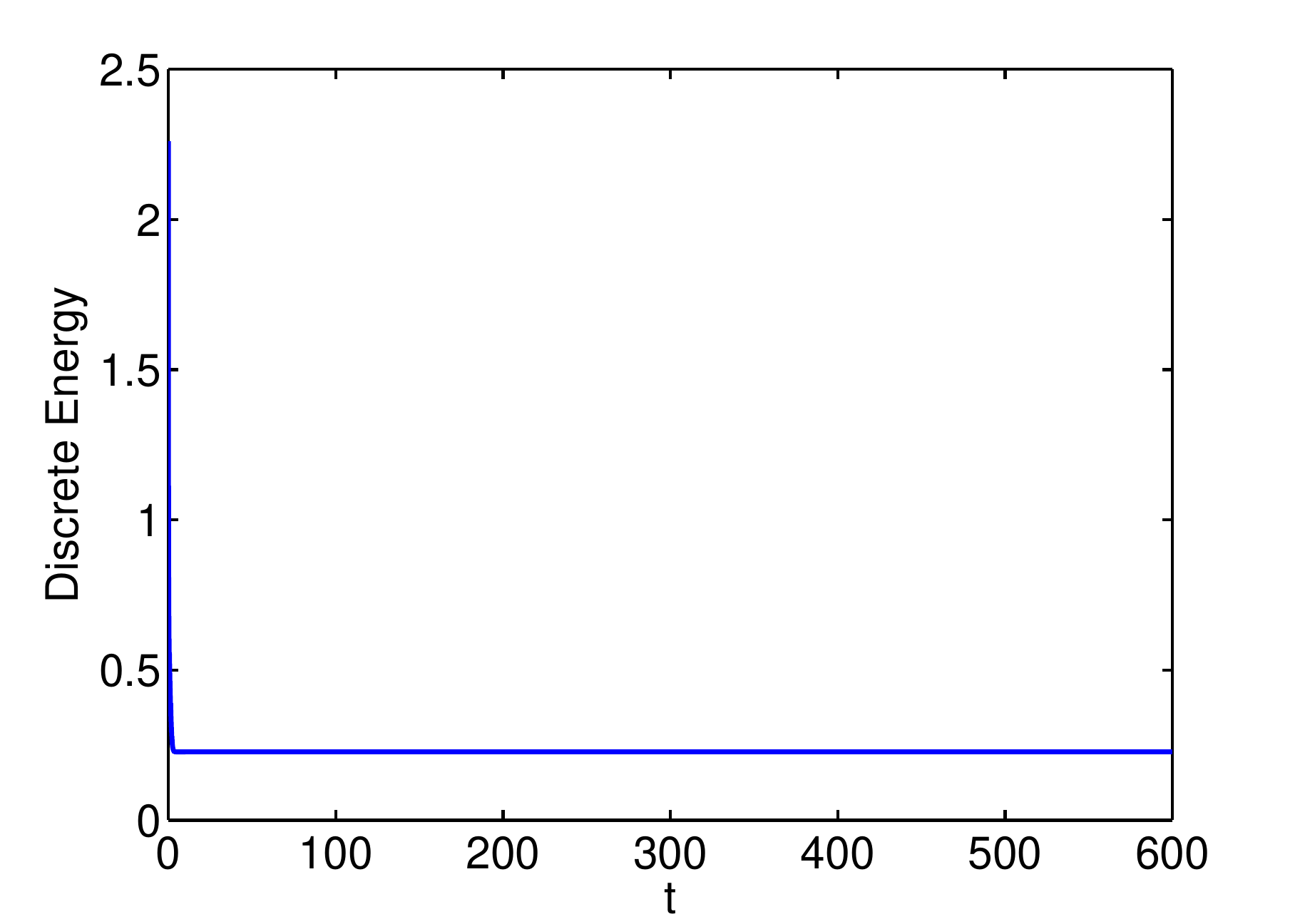}
	\includegraphics[height=0.4\textwidth, width=0.4\textwidth] {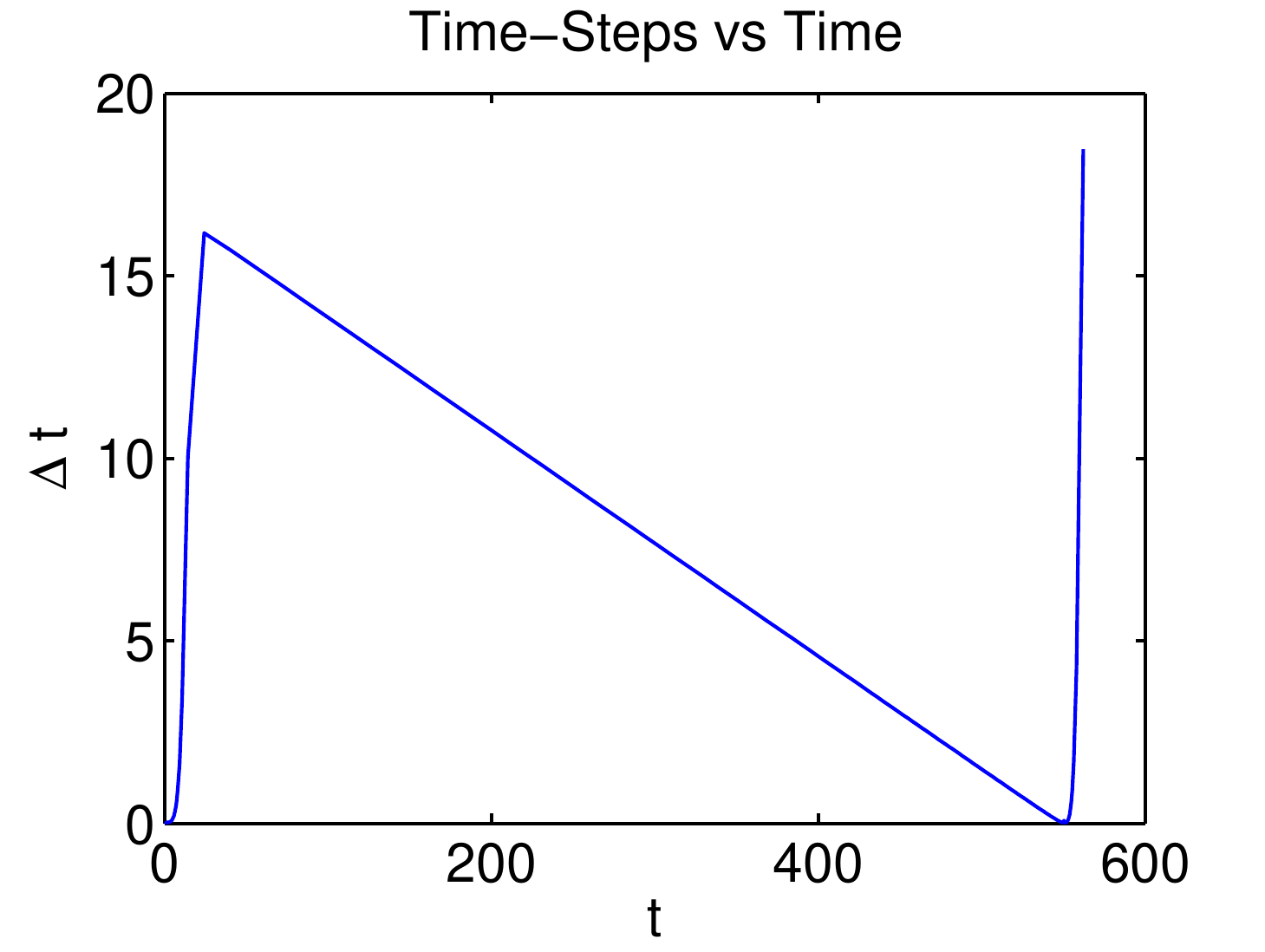}
	\caption{Example \ref{ex2}: Energy decrease and evolution of time steps.\label{1DPeriodic_energy}}
\end{figure}

We consider the 1D Allen-Cahn equation \eqref{allencahn} with the initial condition
$u(x,0)  =  0.8 + \sin(x)$, constant mobility  $\mu(u)=1$ and diffusion constant $\epsilon=0.12$  in the domain $(x,t)\in [0,2\pi] \times [0,600]$. We use the spatial mesh size $\Delta x=\pi /50$. The same problem was solved in  \cite{willoughby11hor} using Fourier spectral space discretization  and with adaptive time integration, as well, using the time integrator pair Backward Differential formula (BDF3) and Adams-Bashforth method (AB-3). For the quartic double-well potential \eqref{func2}, Allen-Cahn equation has one stable ($u=0$) and two unstable ($u=\pm 1$) equilibria, whereas the solutions move from one equilibrium to the other one, which is known as phase separation. The interfaces between two unstable equilibria move over exponentially long times between the region, which is known as the metastability phenomenon. We see in Figure~\ref{1DPeriodic} the fast dynamics from the initial condition to the metastable state, where two transition layers are formed. The numerical  energy is decreasing  in Figure~\ref{1DPeriodic_energy} (left) and time steps are decreased until the metastable state is formed around $t= 546$, Figure~\ref{1DPeriodic_energy} (right).


For the ripening time estimates, the number of time steps are expected to increase linearly by  decreasing the tolerance $\delta_{TOL}$  with the ratio
$\frac{M(\delta_{TOL}/10)}{M(\delta_{TOL})} = \sqrt[p+1]{10}$ \cite{willoughby11hor}.
Table \ref{ex2} shows the ripening time for linear and quadratic DG polynomials. We see that the solution converges at time $T_r=546.5$ with the converge ratio around $\sqrt{10}$.

\begin{table}[htb!]
\caption{Example \ref{ex2}: Convergence of the ripening time with the adaptive AVF method using linear (quadratic) polynomials.\label{table_ex2}}
\centering
\begin{tabular}{c c c c}
$\delta_{TOL}$ & Ripening Time & \# Time Steps & $M(\delta_{TOL}^n)/M(\delta_{TOL}^{n-1})$ \\ \hline
1e-04 &  549.52 (539.71) &  480 (480) &3.02 (3.02)\\
1e-05 & 554.46 (544.54) &  1515 (1515)& 3.12 (3.16)\\
1e-06 &   555.99 (546.05) &  4792 (4790) &3.16 (3.16)\\
1e-07 & 556.47 (546.52)  &  15153 (15152)&3.16 (3.16) \\ \hline
\end{tabular}
\end{table}

\subsection{2D Allen--Cahn equation with constant mobility and double-well potential}\label{ex3}
\begin{figure}[htb!]
\centering
\includegraphics[height=0.3\textwidth, width=0.4\textwidth] {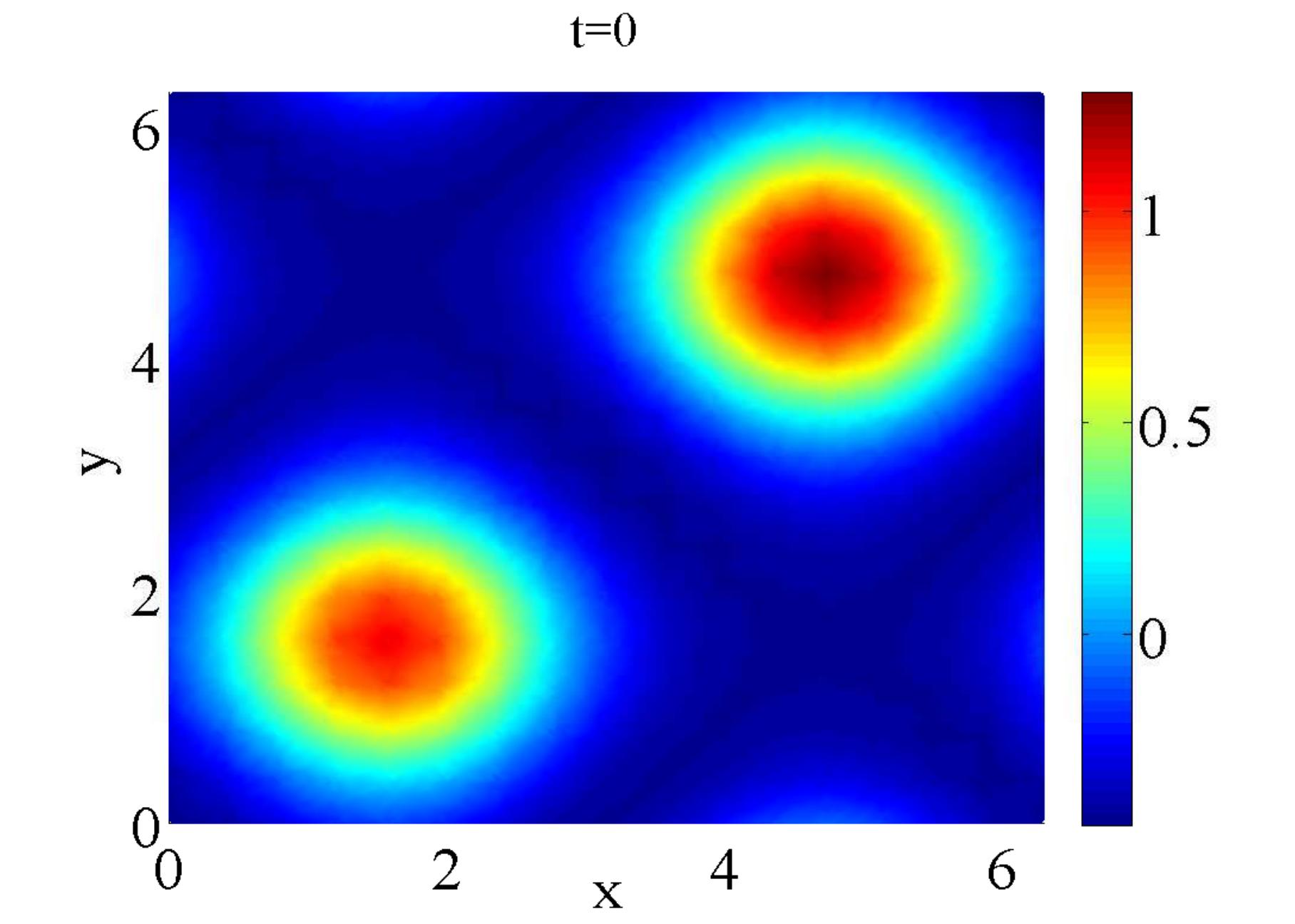}
\includegraphics[height=0.3\textwidth, width=0.4\textwidth] {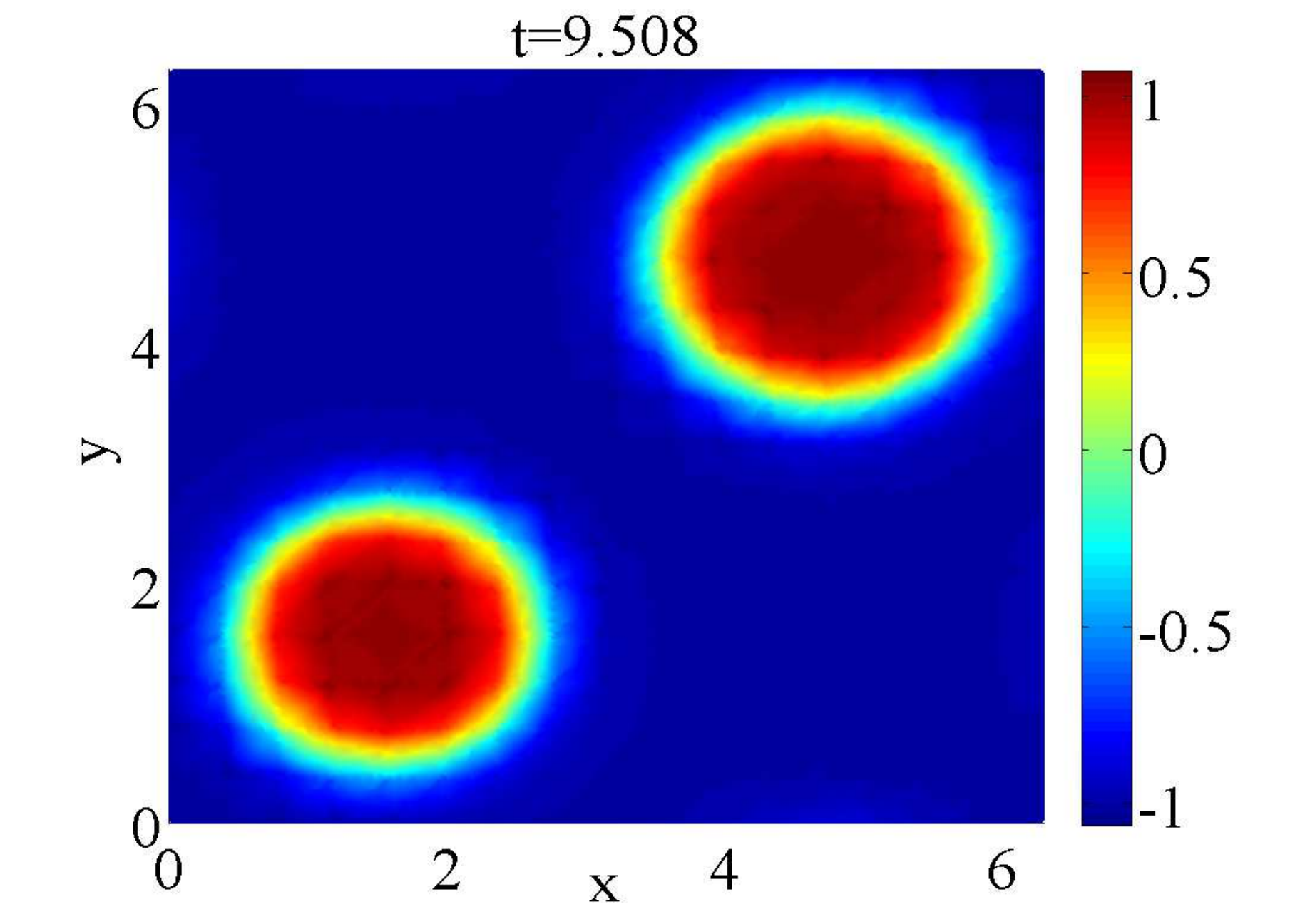}
\includegraphics[height=0.3\textwidth, width=0.4\textwidth] {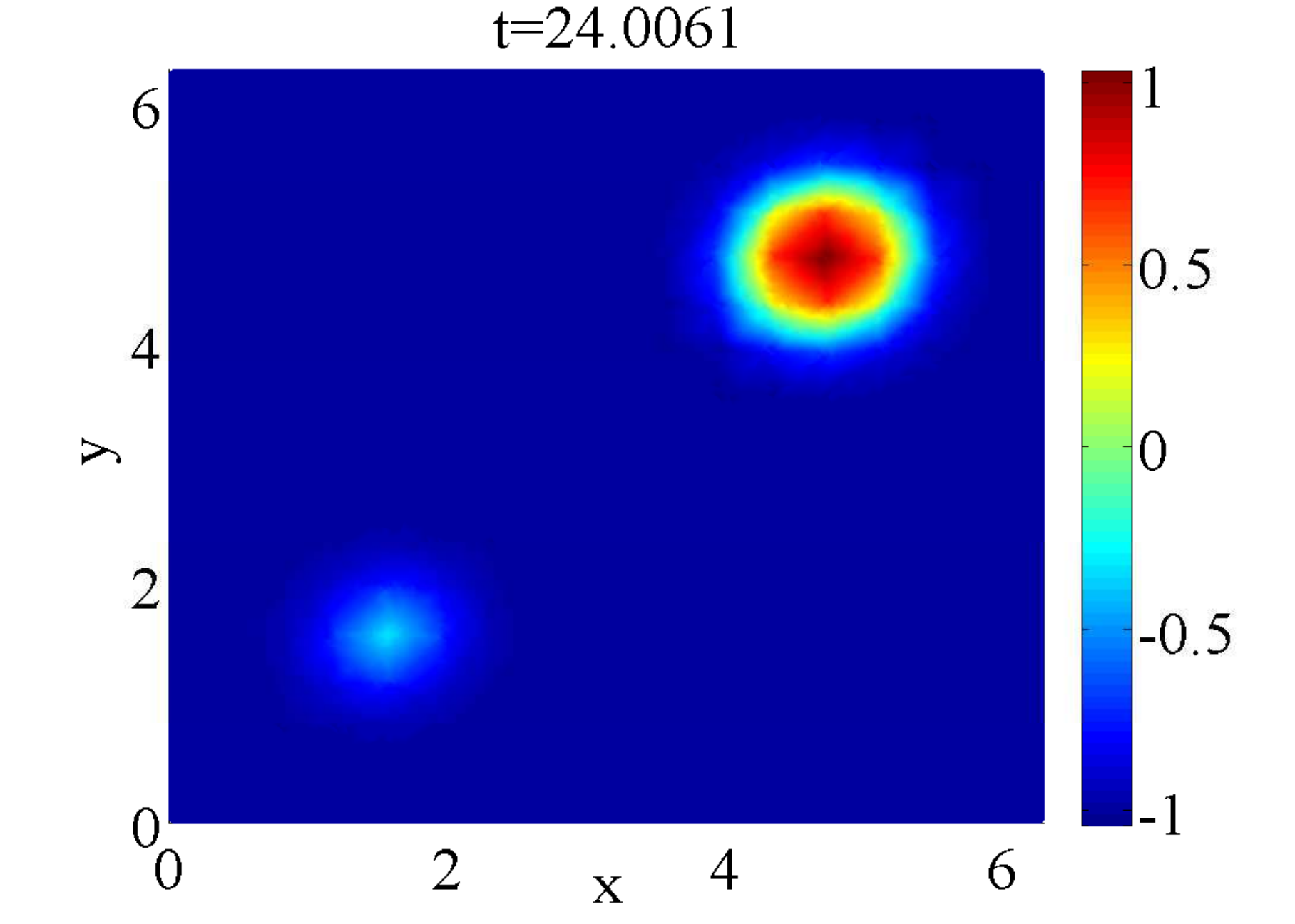}
\includegraphics[height=0.3\textwidth, width=0.4\textwidth] {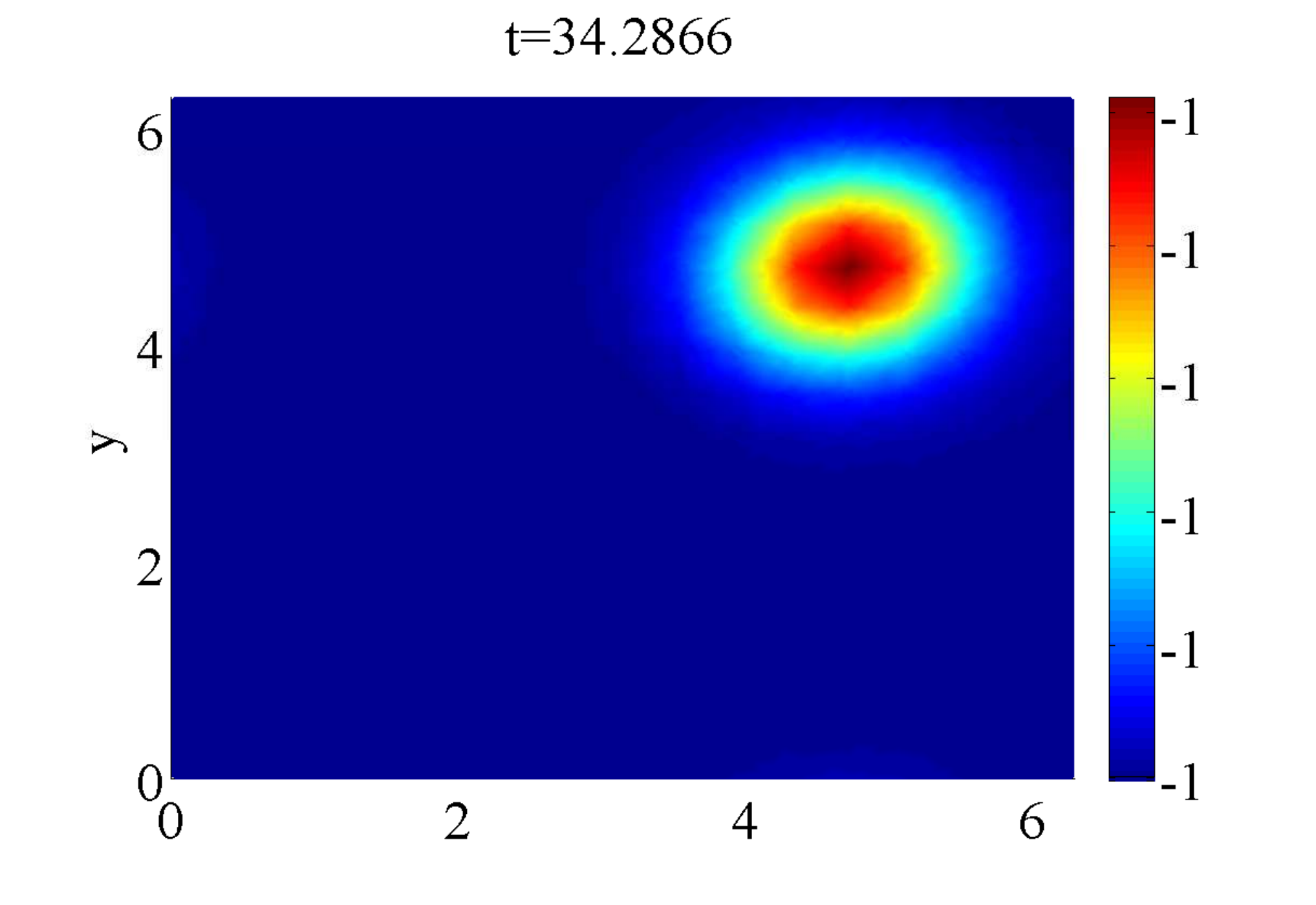}
\caption{Example \ref{ex3}: Evolutions of solution.\label{2Dperiodic_const_mob}}
\end{figure}

We consider the 2D Allen--Cahn equation \eqref{allencahn} with the initial condition \cite{christlieb14has,willoughby11hor}
$$
u(x,y,0) = 2e^{\sin (x) + \sin(y) - 2} + 2.2e^{-\sin (x) - \sin(y) - 2} + 1,
$$
constant mobility  $\mu (u)=1$ and the diffusion constant $\epsilon=0.18$  in the domain $(x,y,t)\in  [0,2\pi]^2 \times [0,33] $. We take as the mesh size $\Delta x = \Delta y=\pi /8$. The solutions with contour plots obtained by the time adaptive scheme with the initial time step size $\tau =0.05$  are shown in Figure~\ref{2Dperiodic_const_mob}. The smaller region is annihilated prior to the larger region. Both regions reach the stable state of $u=-1$ at the end as we expect. We clearly see that the circle shrinks as theoretically predicted, which agrees with numerical results in  \cite{christlieb14has}.
\begin{figure}[htb!]
\centering
\includegraphics[height=0.3\textwidth, width=0.4\textwidth] {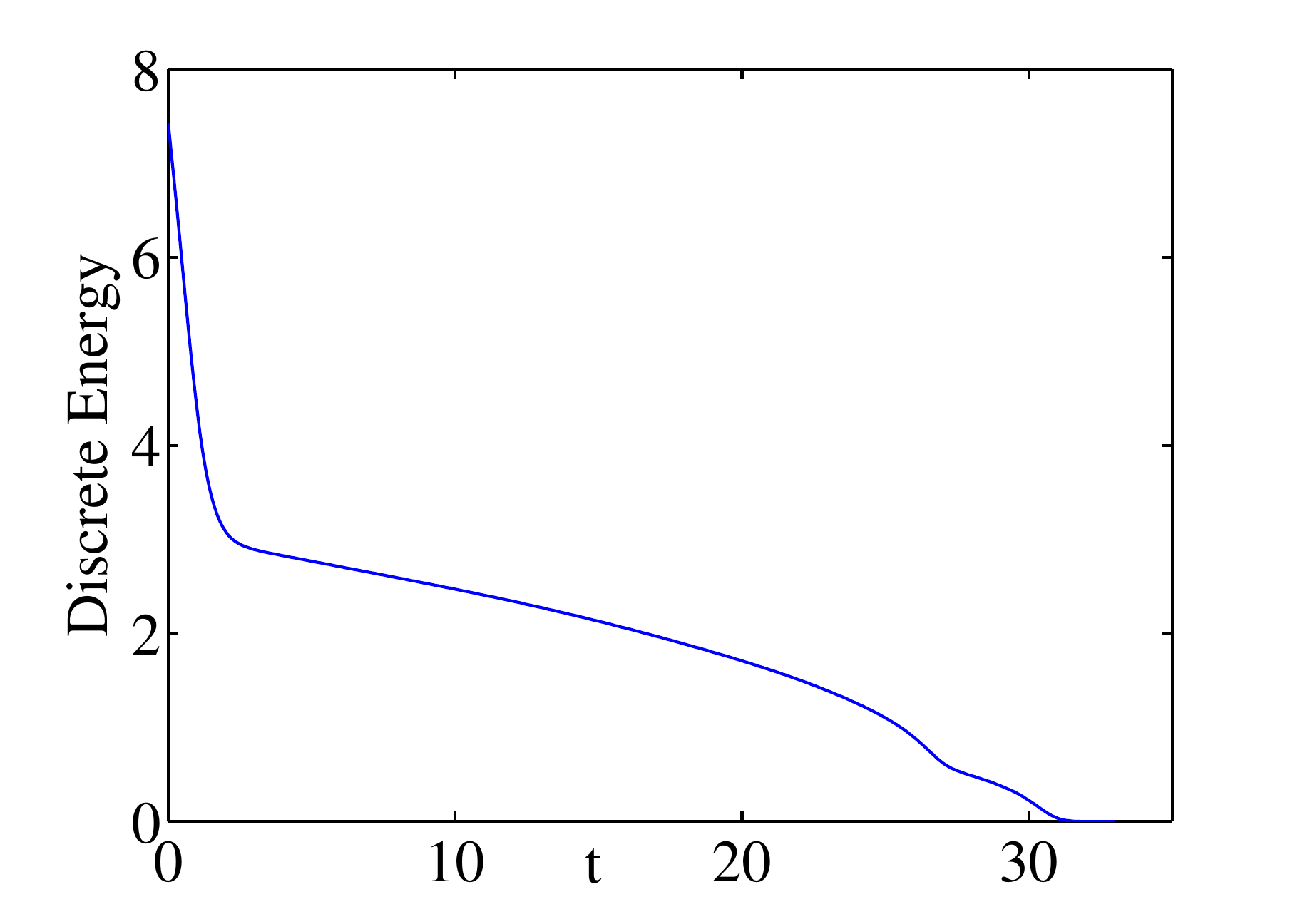}
\includegraphics[height=0.3\textwidth, width=0.4\textwidth] {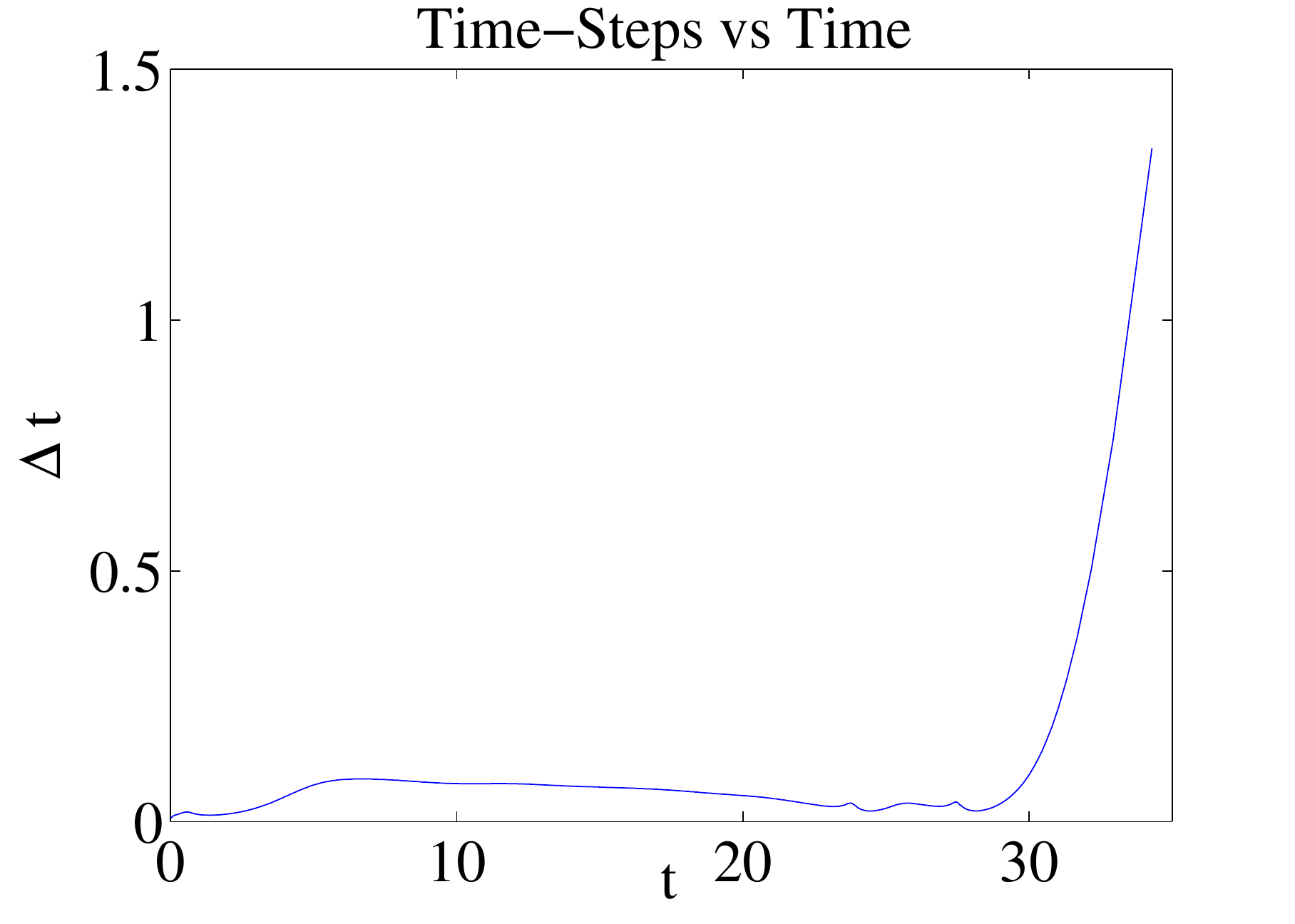}
\caption{Example \ref{ex3}: Energy decrease and time step change. \label{2Dperiodic_const_mobet}}	
\end{figure}

The ripening times for different tolerances with linear and quadratic polynomials are given in Table \ref{table_ex4}. We observe that the ripening time converges by decreasing tolerance and the ratio is close to the theoretically expected value $\sqrt{10}$. The numerical energy is also decreasing monotonically in Figure \ref{2Dperiodic_const_mobet} (left). Small time steps are required until formation of the metastable state around $t=30$, afterward time steps are increased, Figure~\ref{2Dperiodic_const_mobet} (right). Similar results are obtained for the Allen-Cahn and Cahn-Hilliard equations  in \cite{christlieb14has,willoughby11hor}.

\begin{table}[htb!]
\caption{Example \ref{ex3}: Convergence of the ripening time with adaptive AVF method using linear (quadratic) polynomials.\label{table_ex4}}
\centering
\begin{tabular}{c c c c}
$\delta_{TOL}$ & Ripening Time & \# Time Steps  & $M(\delta_{TOL}^n)/M(\delta_{TOL}^{n-1})$  \\ \hline
1e-03 &   27.20 (30.10) &  209 (216)& 3.12 (3.13)\\
1e-04 & 27.33 (30.24) &   668  (692) & 3.20 (3.20)\\
1e-05 &   27.37 (30.25) &  2121  (2197)&3.18 (3.17)\\
1e-06 & 27.37 (30.27)  &  6707 (6956)&3.16 (3.17) \\ \hline
\end{tabular}
\end{table}

\subsection{2D Allen--Cahn equation with constant mobility and logarithmic free energy}\label{ex5}
We consider, as in \cite{shen14omp}, the 2D Allen-Cahn equation \eqref{allencahn} with constant mobility $\mu(u)=2$ and  the diffusion constant $\epsilon=0.04$ in the domain $\Omega=[0,2\pi]\times[0,2\pi]$ for $t\in [0,10]$. The initial condition is $u_0(x,0) = 0.05(2\times \hbox{rand} - 1)$ where 'rand' stands for a random numbers in $[0,1]$.

\begin{figure}[htb!]
\centering
\includegraphics[height=0.3\textwidth, width=0.4\textwidth]{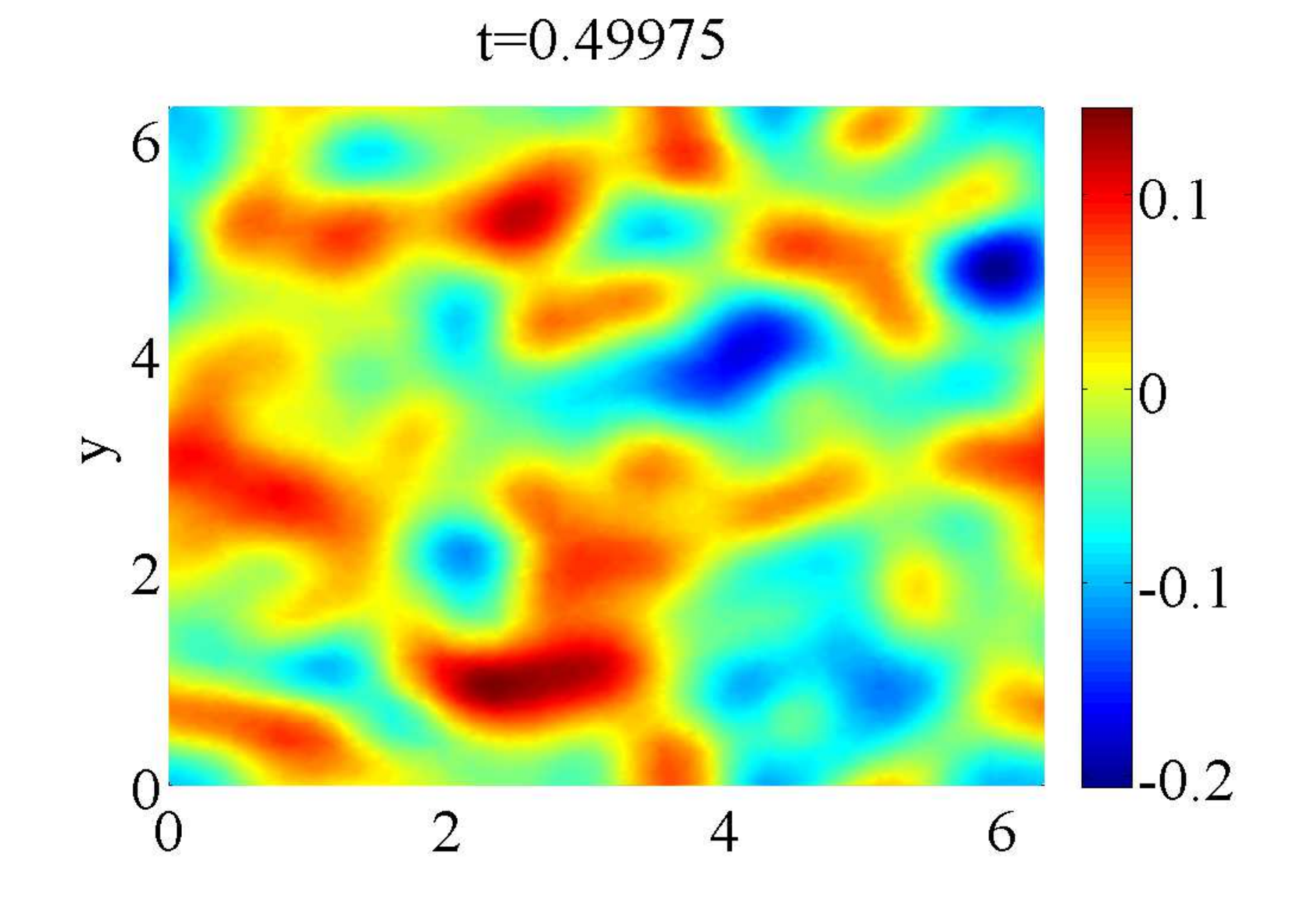}
\includegraphics[height=0.3\textwidth, width=0.4\textwidth]{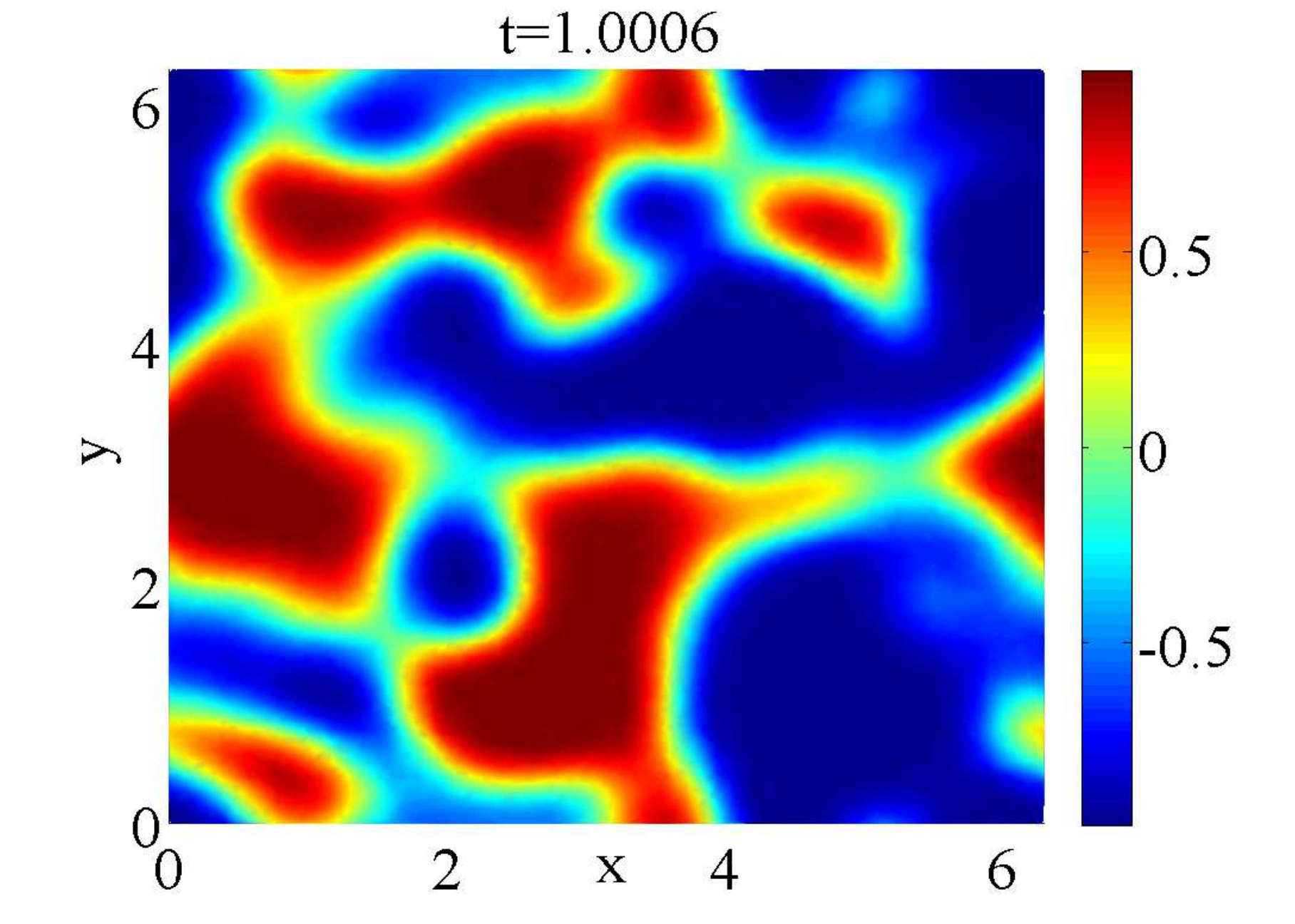}
\includegraphics[height=0.3\textwidth, width=0.4\textwidth]{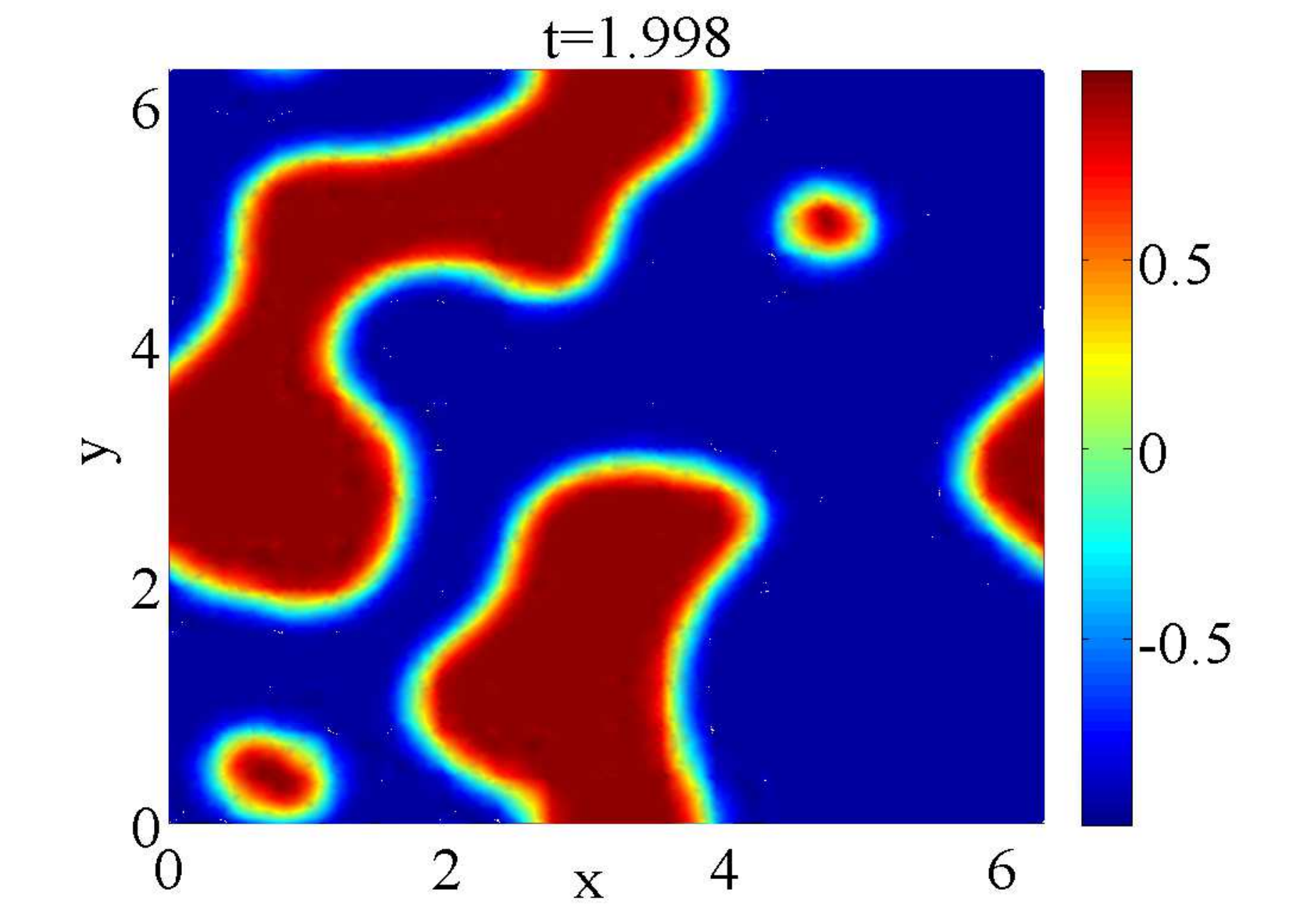}
\includegraphics[height=0.3\textwidth, width=0.4\textwidth]{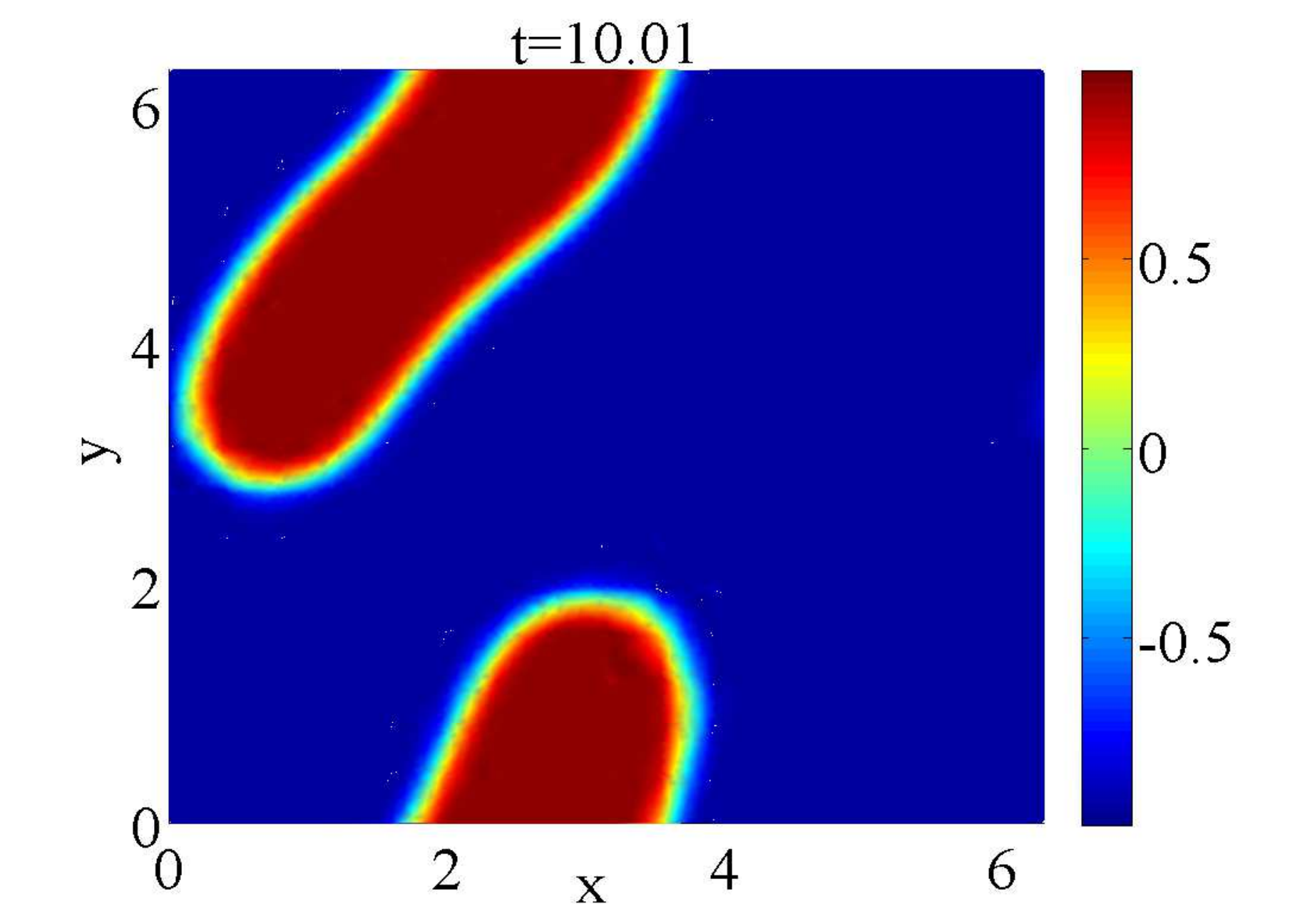}
\caption{Example \ref{ex5}: Evolution of solutions.\label{2D_Periodic_Solution_log}}
\end{figure}

\begin{figure}[htb!]
\centering
\includegraphics[height=0.3\textwidth, width=0.4\textwidth] {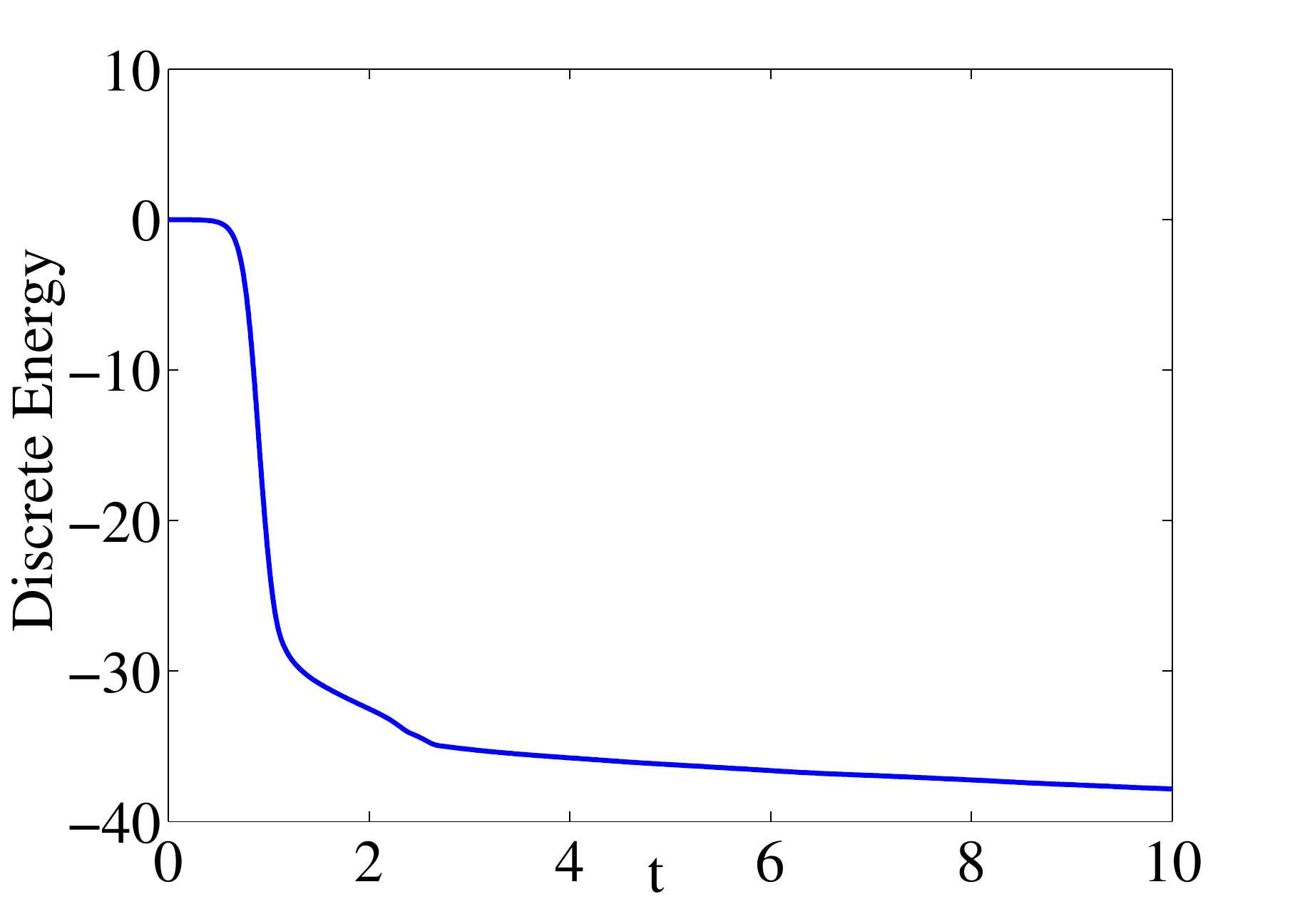}
\includegraphics[height=0.3\textwidth, width=0.4\textwidth] {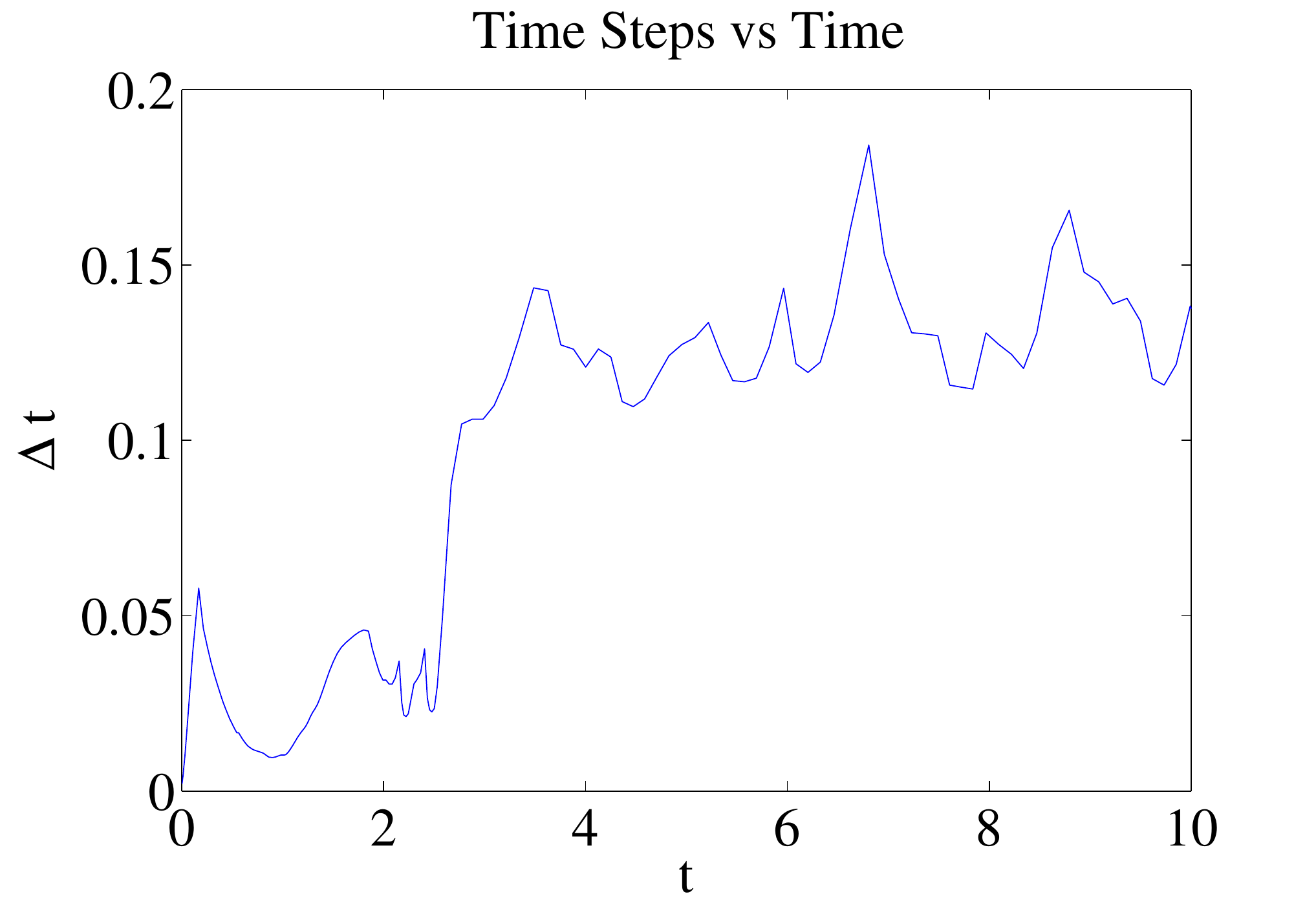}
\caption{Example \ref{ex5}: Energy decrease and time step change.\label{2D_Periodic_Solution_loget}}
\end{figure}

The snapshots of phase evolution is obtained for parameter values $\theta=0.15, \theta_c=0.30$ with time adaptive scheme.
The coarsening phenomena can be seen clearly in Figure~\ref{2D_Periodic_Solution_log}. The numerical energy decrease by the time is seen clearly in Figure~\ref{2D_Periodic_Solution_loget}.


\subsection{2D Allen--Cahn with degenerate mobility and logarithmic free energy}\label{ex7}
We consider the $2D$ Allen--Cahn equation \eqref{allencahn} with the degenerate mobility  $\mu(u)=2(1-u^2)$  \cite{shen14omp} and the diffusion constant $\epsilon=0.04$ in the domain $\Omega=[0,2\pi]\times[0,2\pi]$ for $t\in [0,10]$.  The initial condition is $u_0(x,0) = 0.05(2\times \hbox{rand} - 1)$.

\begin{figure}[htb!]
\centering
\includegraphics[height=0.3\textwidth, width=0.4\textwidth] {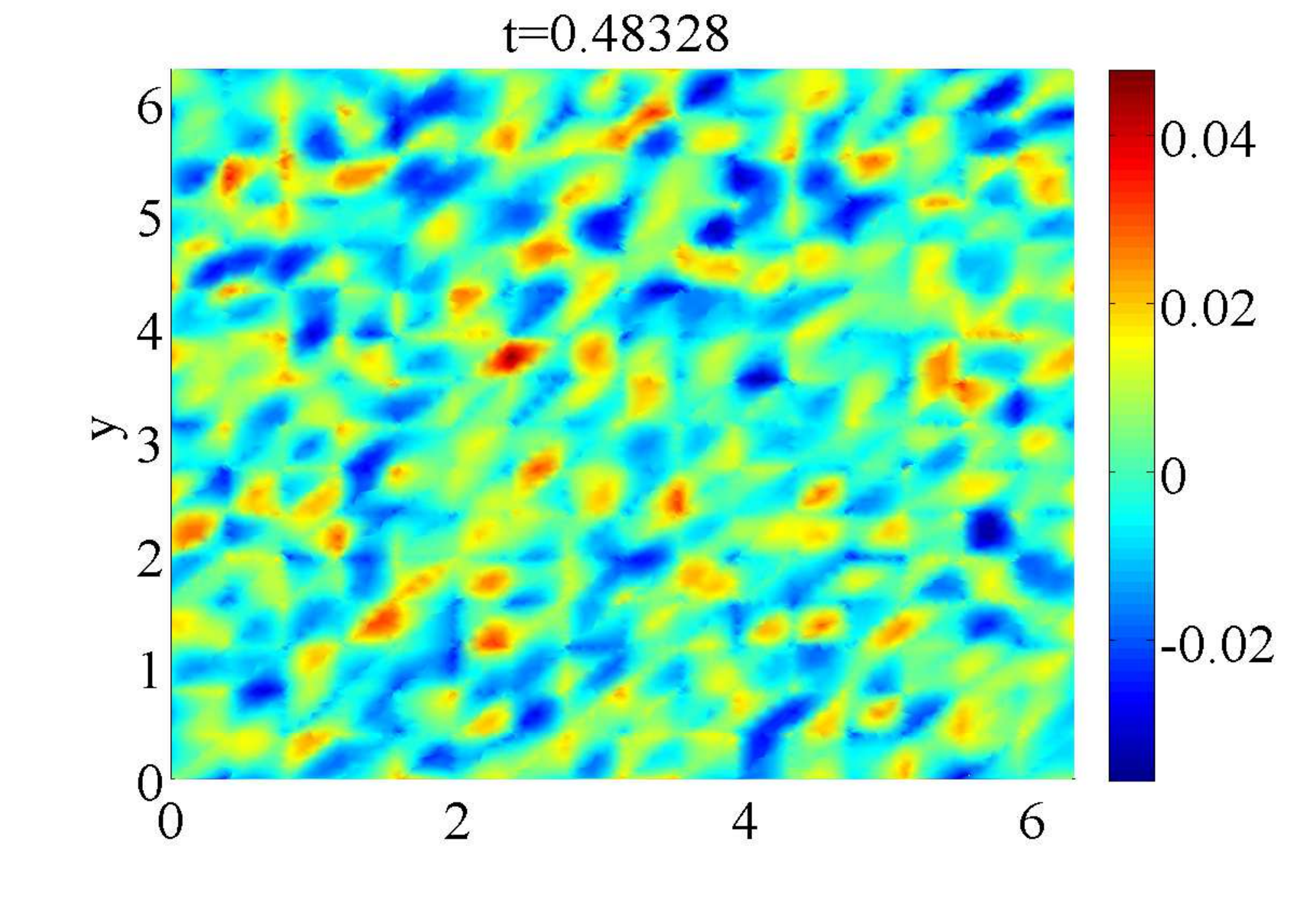}
\includegraphics[height=0.3\textwidth, width=0.4\textwidth]{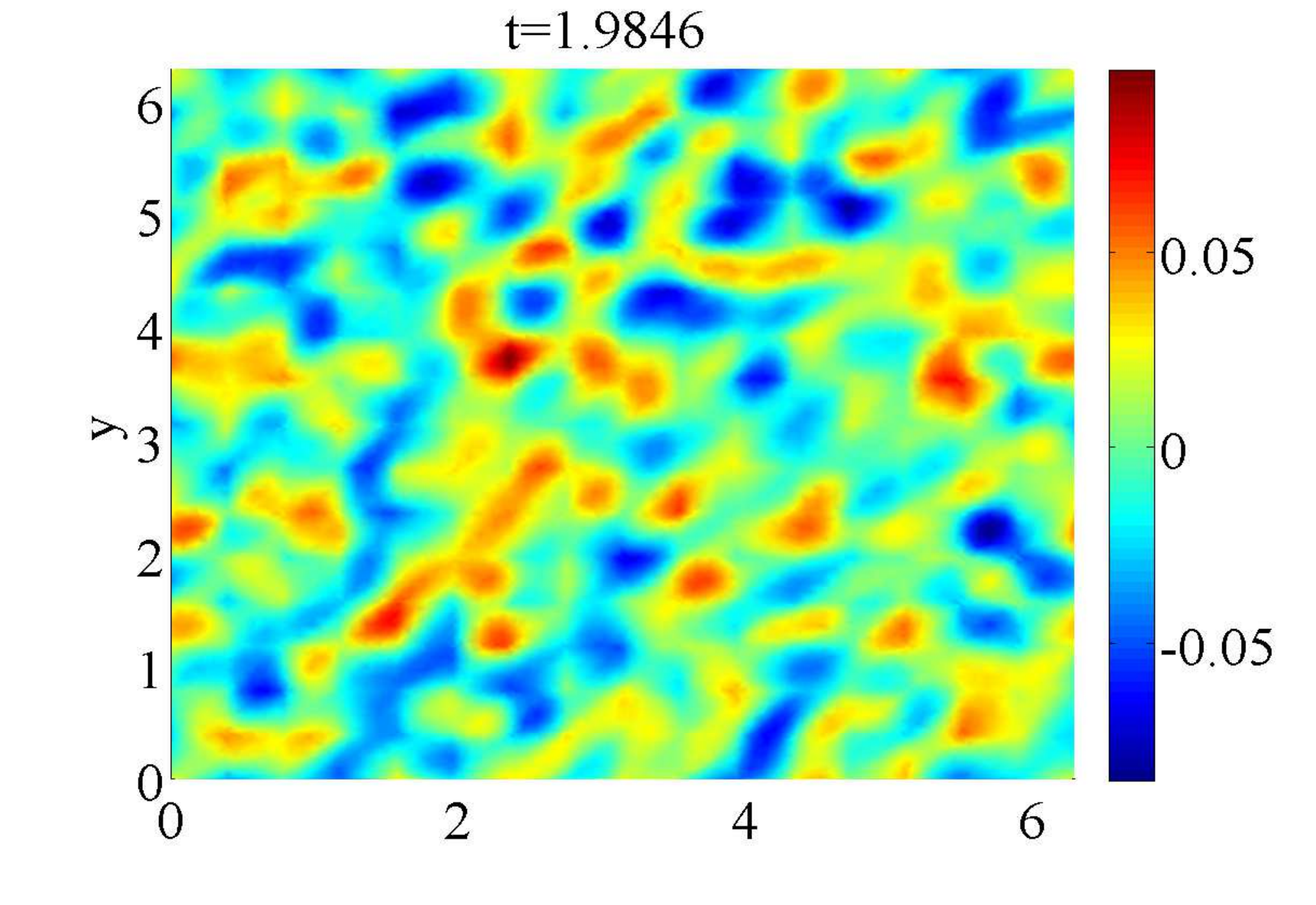}
\includegraphics[height=0.3\textwidth, width=0.4\textwidth]{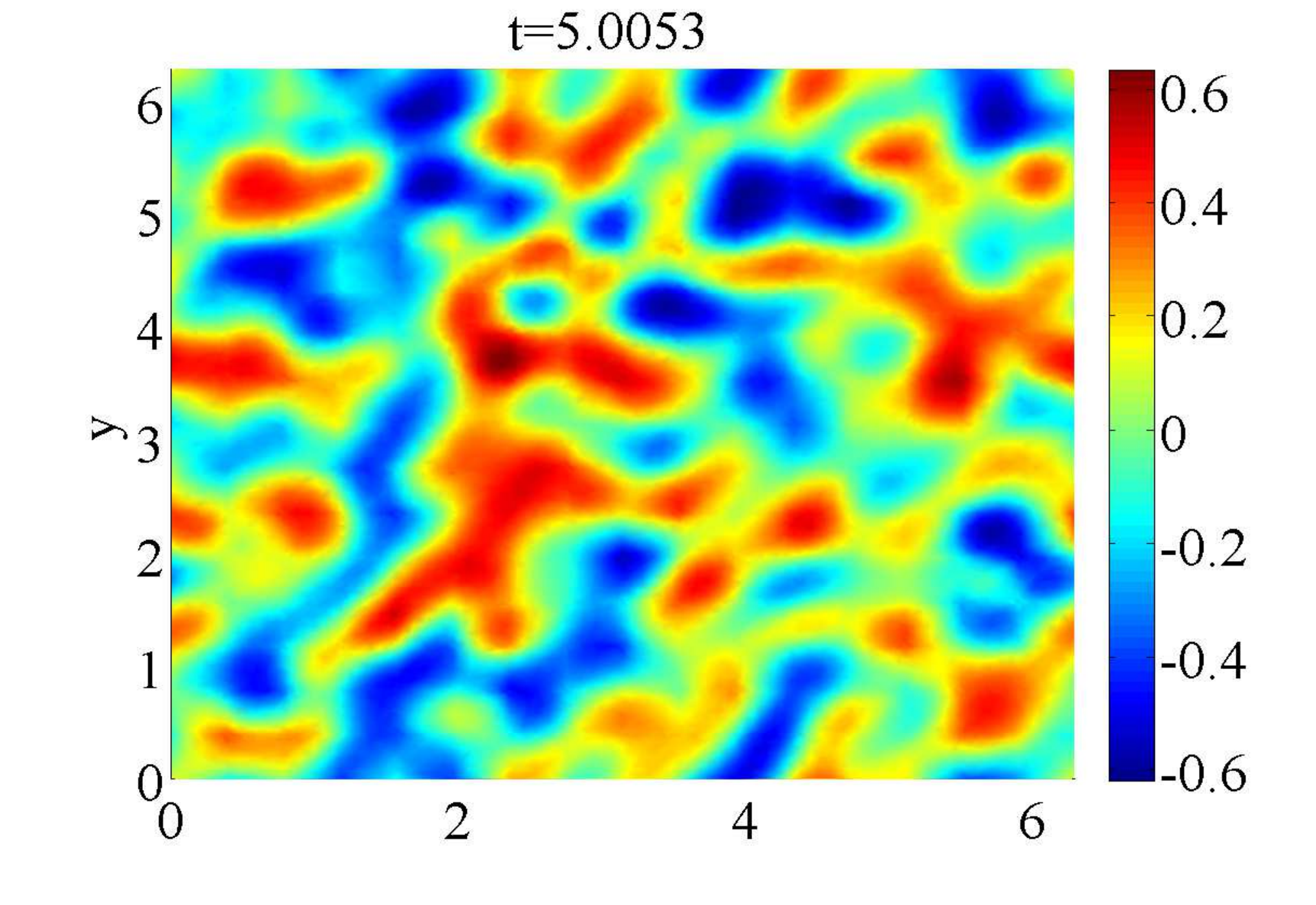}
\includegraphics[height=0.3\textwidth, width=0.4\textwidth]{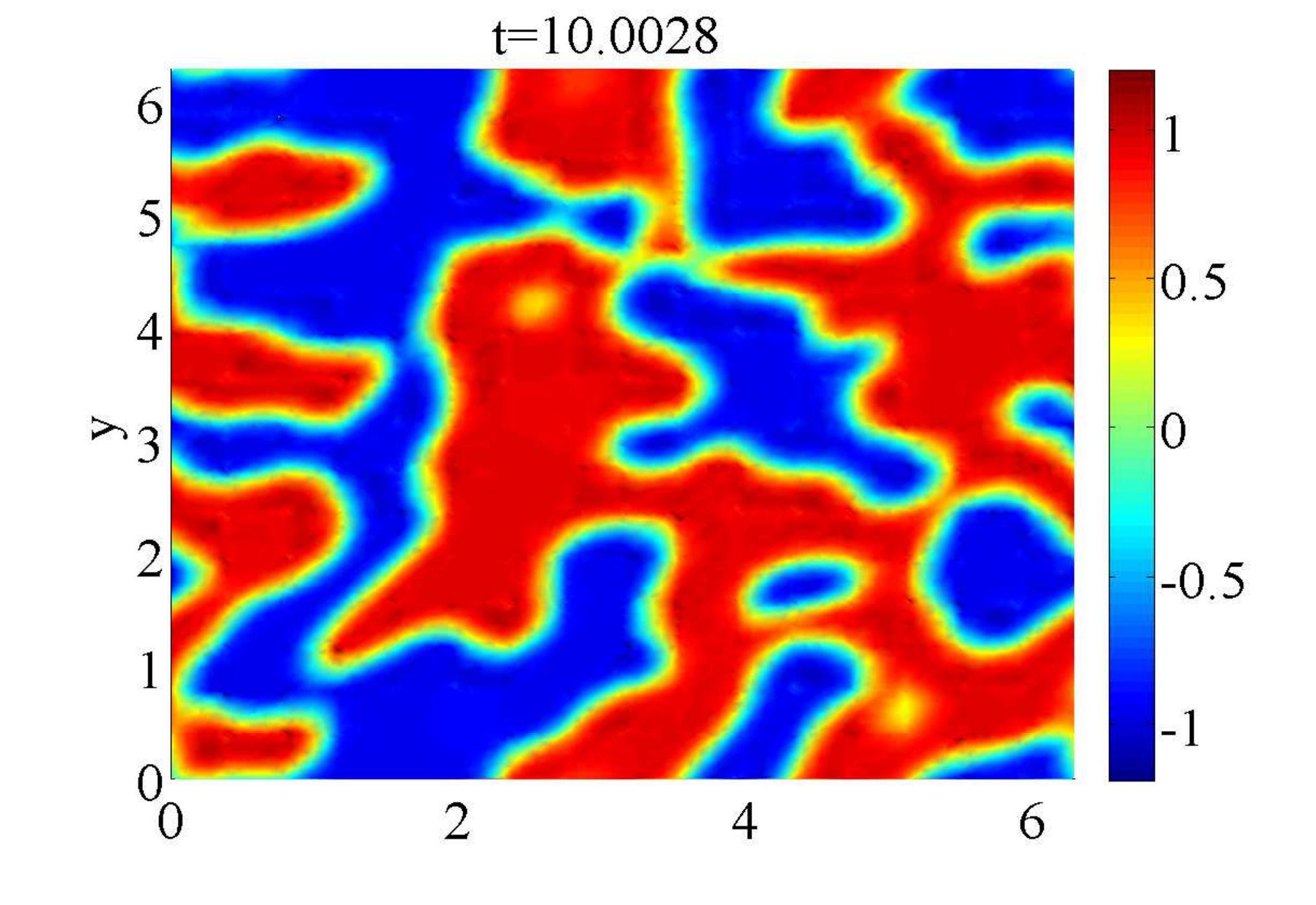}
\caption{Example \ref{ex7}: Evolutions of solution. \label{2D_Periodic_Solution_log_mobility_adap}}
\end{figure}

\begin{figure}[htb!]
\centering
\includegraphics[height=0.3\textwidth, width=0.4\textwidth]{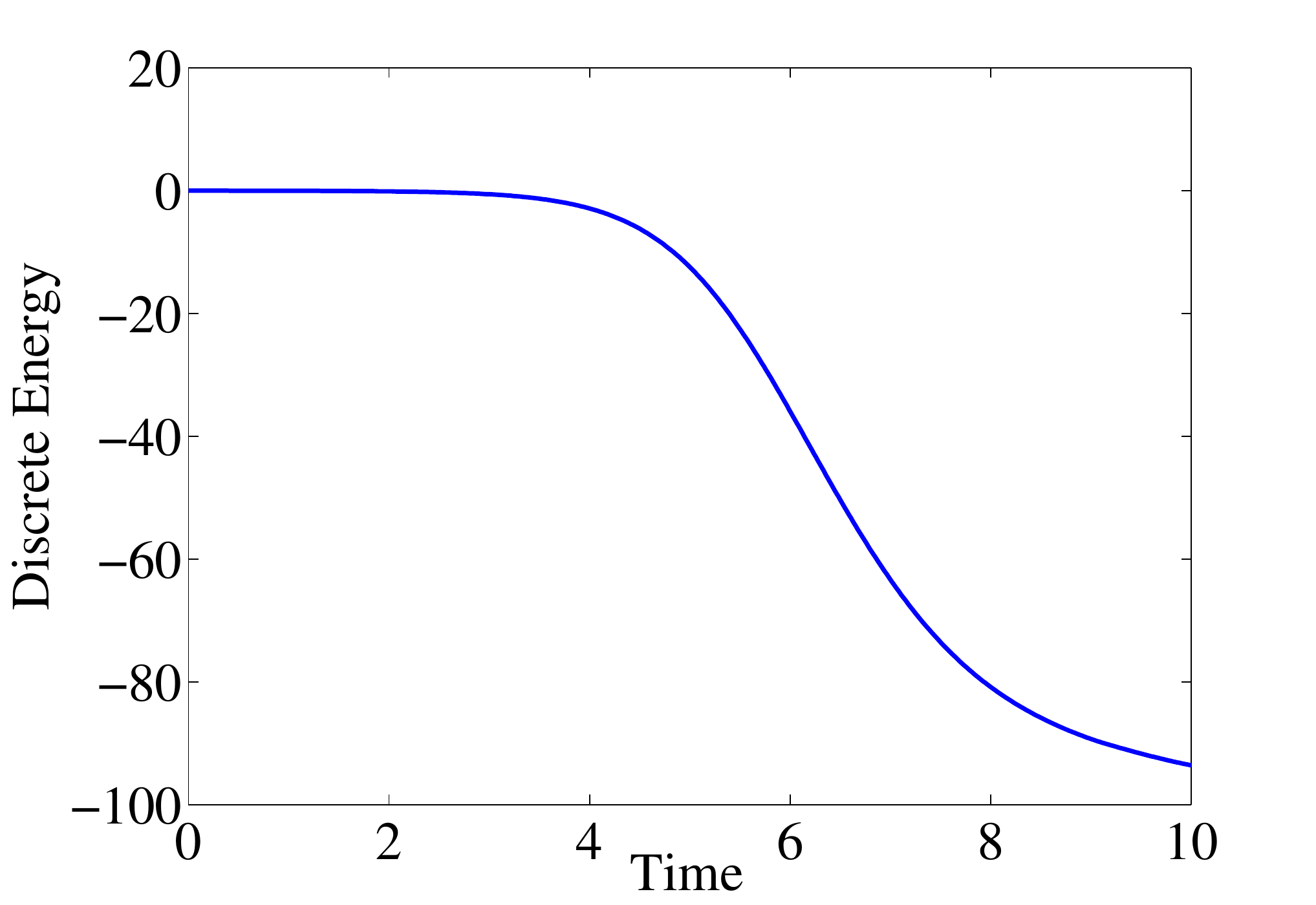}
\includegraphics[height=0.3\textwidth, width=0.4\textwidth]{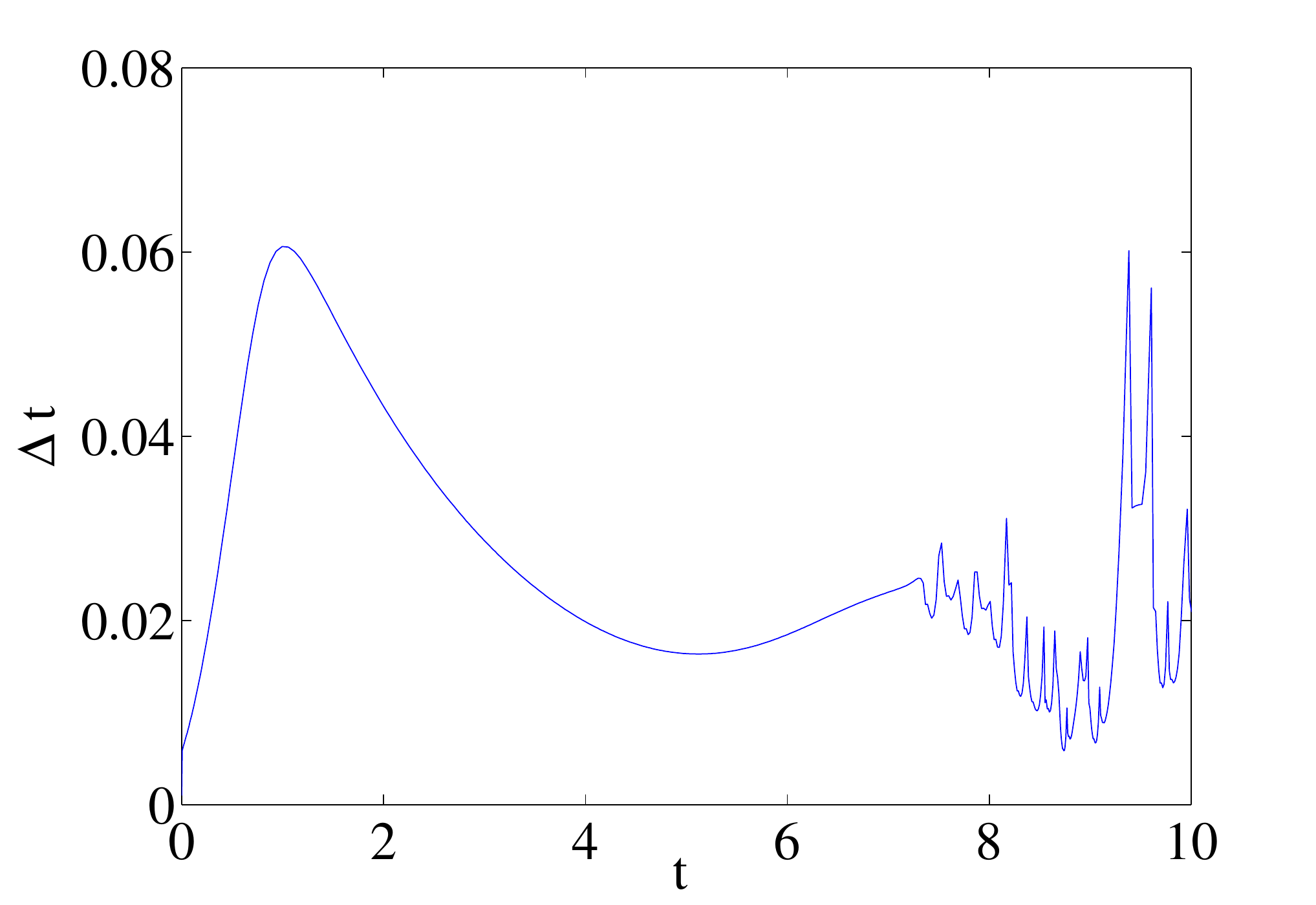}
\caption{Example \ref{ex7}: Energy decrease and time step change.\label{2D_Periodic_Solution_log_mobility_adapet}}
\end{figure}

The phase evolution is obtained for parameter values $\theta=0.50, \theta_c=0.95$. In Figure~\ref{2D_Periodic_Solution_log_mobility_adap}, the corresponding solution contours are plotted, while the numerical energy decrease can be seen in Figure~\ref{2D_Periodic_Solution_log_mobility_adapet}. The numerical results are similar to those in \cite{shen14omp}.


\section{Conclusions}\label{section5}
Numerical results for one and two dimensional Allen-Cahn equation with constant and degenerate mobility, and with polynomial and logarithmic free energy  illustrate the applicability of the SIPG and AVF methods with adaptive time stepping to resolve accurately the dynamics of the Allen-Cahn equation. The ripening time can be detected correctly and the metastability phenomena can be observed numerically. Because the DG method is suitable for handling sharp interfaces and singularities due to its local nature, in a future work, we will study the adaptive DGFEM methods for the sharp interface limit ($\epsilon \rightarrow 0)$ of the Allen-Cahn equation.

\section*{Acknowledgments}
This work has been supported by Scientific HR Development Program (\"OYP) of the Turkish Higher Education Council (Y\"OK).


\begin{thebibliography}{10}

\bibitem{allen79amt}
M.~S. Allen and J.~W. Cahn.
\newblock A microscopic theory for antiphase boundary motion and its
  application to antiphase domain coarsening.
\newblock {\em Acta Metallurgica}, 27(6):1085--1095, 1979.

\bibitem{arnold82ipf}
D.~N. Arnold.
\newblock An interior penalty finite element method with discontinuous
  elements.
\newblock {\em SIAM J. Numer. Anal.}, 19:724--760, 1982.

\bibitem{barrett99fea}
J.~W. Barrett and J.~F. Blowey.
\newblock Finite element approximation of the {C}ahn--{H}illiard equation with
  concentration dependent mobility.
\newblock {\em Mathematics of Computation}, 68(226):487--517, 1999.

\bibitem{barrett00fea}
J.~W. Barrett, J.~F. Blowey, and H.~Garcke.
\newblock Finite element approximation of the {C}ahn--{H}illiard equation with
  degenerate mobility.
\newblock {\em SIAM Journal on Numerical Analysis}, 37(1):286--318, 2000.

\bibitem{celledoni12ped}
E.~Celledoni, V.~Grimm, R.I. McLachlan, D.I. McLaren, D.~O’Neale, B.~Owren,
  and G.R.W. Quispel.
\newblock Preserving energy resp. dissipation in numerical {PDEs} using the
  “average vector field” method.
\newblock {\em Journal of Computational Physics}, 231(20):6770 -- 6789, 2012.

\bibitem{christlieb14has}
A.~Christlieb, J.~Jones, B.~Wetton K.~Promislow, and M.~Willoughby.
\newblock High accuracy solutions to energy gradient flows from material
  science models.
\newblock {\em Journal of Computational Physics}, 257:193--215, 2014.

\bibitem{deuflhard12ans}
P.~Deuflhard and M.~Weisser.
\newblock {\em Adaptive Numerical Solutions of Partial Differential Equations}.
\newblock de Gruyter, Berlin, 2012.

\bibitem{feng14aip}
X.~Feng and Y.~Li.
\newblock Analysis of interior penalty discontinuous {G}alerkin methods for the
  {A}llen--{C}ahn equation and the mean curvature flow.
\newblock {\em arXiv:1310.7504v2 [math.NA}, 2014.

\bibitem{feng13nsi}
X.~Feng, H.~Song, T.~Tang, and J.~Yang.
\newblock Nonlinear stability of the implicit--explicit methods for the
  {A}llen--{C}ahn equation.
\newblock {\em Inverse Problems and Imaging}, 7(3):679--695, 2013.

\bibitem{feng13scn}
Xinlong Feng, Tao Tang, and Jiang Yang.
\newblock Stabilized {Crank-Nicolson/Adams-Bashforth} schemes for phase field
  models.
\newblock {\em East Asian J. Appl. Math}, 3:59--80, 2013.

\bibitem{Guo14esd}
R.~Guo and Y.~Xu.
\newblock Efficient solvers of discontinuous {G}alerkin discretization for the
  {C}ahn-{H}illiard equations.
\newblock {\em Journal of Scientific Computing}, 58(2):380--408, 2014.

\bibitem{guo-nse}
Ruihan Guo, Liangyue Ji, and Yan Xu.
\newblock Numerical simulation and error estimates for the local discontinuos
  {G}alerkin method of the {A}llen-{C}ahn equation.
\newblock available at http://home.ustc.edu.cn/~guoguo88/, 2013.

\bibitem{hairer10epv}
E.~Hairer.
\newblock Energy-preserving variant of collocation methods.
\newblock {\em Journal of Numerical Analysis, Industrial and Applied
  Mathematics}, 5:73--84, 2010.

\bibitem{hairer96sod}
E.~Hairer and G.~Wanner.
\newblock {\em Solving Ordinary Differential Equations II: Stiff and
  Differential--Algebraic Problems}.
\newblock Springer Series in Computational Mathematics:Berlin. Springer, 1996.

\bibitem{liu13ssi}
Fei Liu and Jie Shen.
\newblock Stabilized semi-implicit spectral deferred correction methods for
  allen--cahn and cahn--hilliard equations.
\newblock {\em Mathematical Methods in the Applied Sciences}, 2013.

\bibitem{riviere08dgm}
B.~Rivi\`{e}re.
\newblock {\em Discontinuous {G}alerkin methods for solving elliptic and
  parabolic equations, Theory and implementation}.
\newblock SIAM, 2008.

\bibitem{schwieweck10sdg}
F.~Schieweck.
\newblock A stable discontinuous {G}alerkin--{P}etrov time discretization of
  higher order.
\newblock {\em Journal of Numerical Mathematics}, 18:25--27, 2010.

\bibitem{shen14omp}
Jie Shen, Tao Tang, and Jiang Yang.
\newblock On the maximum principle preserving schemes for the generalized
  {A}llen-{C}ahn equation.
\newblock Preprint, 2014.

\bibitem{shen10nac}
Jie Shen and Xiaofeng Yang.
\newblock Numerical approximations of {A}llen-{C}ahn and {C}ahn-{H}illiard
  equations.
\newblock {\em Discrete Continuous Dynamical. Syst. A}, 28:1669--1691, 2010.

\bibitem{zee11goe}
K.G. van~der Zee, J.~Tinsley Oden, S.~Prudhomme, and A.~Hawkins-Daarud.
\newblock Goal-oriented error estimation for {C}ahn–-{H}illiard models of
  binary phase transition.
\newblock {\em Numerical Methods for Partial Differential Equations},
  27(1):160--196, 2011.

\bibitem{vemaganti07dgm}
K.~Vemaganti.
\newblock Discontinuous {G}alerkin methods for periodic boundary value
  problems.
\newblock {\em Numerical Methods for Partial Differential Equations},
  23(3):587--596, 2007.

\bibitem{willoughby11hor}
M.R. Willoughby.
\newblock High-order time--adaptive numerical methods for the {A}llen--{C}ahn
  and {C}ahn--{H}illiard equations.
\newblock Master's thesis, The University of British Columbia, Vancouver,
  December 2011.

\end{thebibliography}

\end{document}